
\catcode`\@=11


\message{Loading jyTeX fonts...}



\font\vptrm=cmr5 \font\vptmit=cmmi5 \font\vptsy=cmsy5 \font\vptbf=cmbx5

\skewchar\vptmit='177 \skewchar\vptsy='60 \fontdimen16
\vptsy=\the\fontdimen17 \vptsy

\def\vpt{\ifmmode\err@badsizechange\else
     \@mathfontinit
     \textfont0=\vptrm  \scriptfont0=\vptrm  \scriptscriptfont0=\vptrm
     \textfont1=\vptmit \scriptfont1=\vptmit \scriptscriptfont1=\vptmit
     \textfont2=\vptsy  \scriptfont2=\vptsy  \scriptscriptfont2=\vptsy
     \textfont3=\xptex  \scriptfont3=\xptex  \scriptscriptfont3=\xptex
     \textfont\bffam=\vptbf
     \scriptfont\bffam=\vptbf
     \scriptscriptfont\bffam=\vptbf
     \@fontstyleinit
     \def\rm{\vptrm\fam=\z@}%
     \def\bf{\vptbf\fam=\bffam}%
     \def\oldstyle{\vptmit\fam=\@ne}%
     \rm\fi}


\font\viptrm=cmr6 \font\viptmit=cmmi6 \font\viptsy=cmsy6
\font\viptbf=cmbx6

\skewchar\viptmit='177 \skewchar\viptsy='60 \fontdimen16
\viptsy=\the\fontdimen17 \viptsy

\def\vipt{\ifmmode\err@badsizechange\else
     \@mathfontinit
     \textfont0=\viptrm  \scriptfont0=\vptrm  \scriptscriptfont0=\vptrm
     \textfont1=\viptmit \scriptfont1=\vptmit \scriptscriptfont1=\vptmit
     \textfont2=\viptsy  \scriptfont2=\vptsy  \scriptscriptfont2=\vptsy
     \textfont3=\xptex   \scriptfont3=\xptex  \scriptscriptfont3=\xptex
     \textfont\bffam=\viptbf
     \scriptfont\bffam=\vptbf
     \scriptscriptfont\bffam=\vptbf
     \@fontstyleinit
     \def\rm{\viptrm\fam=\z@}%
     \def\bf{\viptbf\fam=\bffam}%
     \def\oldstyle{\viptmit\fam=\@ne}%
     \rm\fi}

\font\viiptrm=cmr7 \font\viiptmit=cmmi7 \font\viiptsy=cmsy7
\font\viiptit=cmti7 \font\viiptbf=cmbx7

\skewchar\viiptmit='177 \skewchar\viiptsy='60 \fontdimen16
\viiptsy=\the\fontdimen17 \viiptsy

\def\viipt{\ifmmode\err@badsizechange\else
     \@mathfontinit
     \textfont0=\viiptrm  \scriptfont0=\vptrm  \scriptscriptfont0=\vptrm
     \textfont1=\viiptmit \scriptfont1=\vptmit \scriptscriptfont1=\vptmit
     \textfont2=\viiptsy  \scriptfont2=\vptsy  \scriptscriptfont2=\vptsy
     \textfont3=\xptex    \scriptfont3=\xptex  \scriptscriptfont3=\xptex
     \textfont\itfam=\viiptit
     \scriptfont\itfam=\viiptit
     \scriptscriptfont\itfam=\viiptit
     \textfont\bffam=\viiptbf
     \scriptfont\bffam=\vptbf
     \scriptscriptfont\bffam=\vptbf
     \@fontstyleinit
     \def\rm{\viiptrm\fam=\z@}%
     \def\it{\viiptit\fam=\itfam}%
     \def\bf{\viiptbf\fam=\bffam}%
     \def\oldstyle{\viiptmit\fam=\@ne}%
     \rm\fi}


\font\viiiptrm=cmr8 \font\viiiptmit=cmmi8 \font\viiiptsy=cmsy8
\font\viiiptit=cmti8
\font\viiiptbf=cmbx8

\skewchar\viiiptmit='177 \skewchar\viiiptsy='60 \fontdimen16
\viiiptsy=\the\fontdimen17 \viiiptsy

\def\viiipt{\ifmmode\err@badsizechange\else
     \@mathfontinit
     \textfont0=\viiiptrm  \scriptfont0=\viptrm  \scriptscriptfont0=\vptrm
     \textfont1=\viiiptmit \scriptfont1=\viptmit \scriptscriptfont1=\vptmit
     \textfont2=\viiiptsy  \scriptfont2=\viptsy  \scriptscriptfont2=\vptsy
     \textfont3=\xptex     \scriptfont3=\xptex   \scriptscriptfont3=\xptex
     \textfont\itfam=\viiiptit
     \scriptfont\itfam=\viiptit
     \scriptscriptfont\itfam=\viiptit
     \textfont\bffam=\viiiptbf
     \scriptfont\bffam=\viptbf
     \scriptscriptfont\bffam=\vptbf
     \@fontstyleinit
     \def\rm{\viiiptrm\fam=\z@}%
     \def\it{\viiiptit\fam=\itfam}%
     \def\bf{\viiiptbf\fam=\bffam}%
     \def\oldstyle{\viiiptmit\fam=\@ne}%
     \rm\fi}


\def\getixpt{%
     \font\ixptrm=cmr9
     \font\ixptmit=cmmi9
     \font\ixptsy=cmsy9
     \font\ixptit=cmti9
     \font\ixptbf=cmbx9
     \skewchar\ixptmit='177 \skewchar\ixptsy='60
     \fontdimen16 \ixptsy=\the\fontdimen17 \ixptsy}

\def\ixpt{\ifmmode\err@badsizechange\else
     \@mathfontinit
     \textfont0=\ixptrm  \scriptfont0=\viiptrm  \scriptscriptfont0=\vptrm
     \textfont1=\ixptmit \scriptfont1=\viiptmit \scriptscriptfont1=\vptmit
     \textfont2=\ixptsy  \scriptfont2=\viiptsy  \scriptscriptfont2=\vptsy
     \textfont3=\xptex   \scriptfont3=\xptex    \scriptscriptfont3=\xptex
     \textfont\itfam=\ixptit
     \scriptfont\itfam=\viiptit
     \scriptscriptfont\itfam=\viiptit
     \textfont\bffam=\ixptbf
     \scriptfont\bffam=\viiptbf
     \scriptscriptfont\bffam=\vptbf
     \@fontstyleinit
     \def\rm{\ixptrm\fam=\z@}%
     \def\it{\ixptit\fam=\itfam}%
     \def\bf{\ixptbf\fam=\bffam}%
     \def\oldstyle{\ixptmit\fam=\@ne}%
     \rm\fi}


\font\xptrm=cmr10 \font\xptmit=cmmi10 \font\xptsy=cmsy10
\font\xptex=cmex10 \font\xptit=cmti10 \font\xptsl=cmsl10
\font\xptbf=cmbx10 \font\xpttt=cmtt10 \font\xptss=cmss10
\font\xptsc=cmcsc10 \font\xptbfs=cmb10 \font\xptbmit=cmmib10

\skewchar\xptmit='177 \skewchar\xptbmit='177 \skewchar\xptsy='60
\fontdimen16 \xptsy=\the\fontdimen17 \xptsy

\def\xpt{\ifmmode\err@badsizechange\else
     \@mathfontinit
     \textfont0=\xptrm  \scriptfont0=\viiptrm  \scriptscriptfont0=\vptrm
     \textfont1=\xptmit \scriptfont1=\viiptmit \scriptscriptfont1=\vptmit
     \textfont2=\xptsy  \scriptfont2=\viiptsy  \scriptscriptfont2=\vptsy
     \textfont3=\xptex  \scriptfont3=\xptex    \scriptscriptfont3=\xptex
     \textfont\itfam=\xptit
     \scriptfont\itfam=\viiptit
     \scriptscriptfont\itfam=\viiptit
     \textfont\bffam=\xptbf
     \scriptfont\bffam=\viiptbf
     \scriptscriptfont\bffam=\vptbf
     \textfont\bfsfam=\xptbfs
     \scriptfont\bfsfam=\viiptbf
     \scriptscriptfont\bfsfam=\vptbf
     \textfont\bmitfam=\xptbmit
     \scriptfont\bmitfam=\viiptmit
     \scriptscriptfont\bmitfam=\vptmit
     \@fontstyleinit
     \def\rm{\xptrm\fam=\z@}%
     \def\it{\xptit\fam=\itfam}%
     \def\sl{\xptsl}%
     \def\bf{\xptbf\fam=\bffam}%
     \def\tt{\xpttt}%
     \def\ss{\xptss}%
     \def\sc{\xptsc}%
     \def\bfs{\xptbfs\fam=\bfsfam}%
     \def\bmit{\fam=\bmitfam}%
     \def\oldstyle{\xptmit\fam=\@ne}%
     \rm\fi}


\def\getxipt{%
     \font\xiptrm=cmr10  scaled\magstephalf
     \font\xiptmit=cmmi10 scaled\magstephalf
     \font\xiptsy=cmsy10 scaled\magstephalf
     \font\xiptex=cmex10 scaled\magstephalf
     \font\xiptit=cmti10 scaled\magstephalf
     \font\xiptsl=cmsl10 scaled\magstephalf
     \font\xiptbf=cmbx10 scaled\magstephalf
     \font\xipttt=cmtt10 scaled\magstephalf
     \font\xiptss=cmss10 scaled\magstephalf
     \skewchar\xiptmit='177 \skewchar\xiptsy='60
     \fontdimen16 \xiptsy=\the\fontdimen17 \xiptsy}

\def\xipt{\ifmmode\err@badsizechange\else
     \@mathfontinit
     \textfont0=\xiptrm  \scriptfont0=\viiiptrm  \scriptscriptfont0=\viptrm
     \textfont1=\xiptmit \scriptfont1=\viiiptmit \scriptscriptfont1=\viptmit
     \textfont2=\xiptsy  \scriptfont2=\viiiptsy  \scriptscriptfont2=\viptsy
     \textfont3=\xiptex  \scriptfont3=\xptex     \scriptscriptfont3=\xptex
     \textfont\itfam=\xiptit
     \scriptfont\itfam=\viiiptit
     \scriptscriptfont\itfam=\viiptit
     \textfont\bffam=\xiptbf
     \scriptfont\bffam=\viiiptbf
     \scriptscriptfont\bffam=\viptbf
     \@fontstyleinit
     \def\rm{\xiptrm\fam=\z@}%
     \def\it{\xiptit\fam=\itfam}%
     \def\sl{\xiptsl}%
     \def\bf{\xiptbf\fam=\bffam}%
     \def\tt{\xipttt}%
     \def\ss{\xiptss}%
     \def\oldstyle{\xiptmit\fam=\@ne}%
     \rm\fi}


\font\xiiptrm=cmr12 \font\xiiptmit=cmmi12 \font\xiiptsy=cmsy10
scaled\magstep1 \font\xiiptex=cmex10  scaled\magstep1
\font\xiiptit=cmti12 \font\xiiptsl=cmsl12 \font\xiiptbf=cmbx12
\font\xiiptss=cmss12 \font\xiiptsc=cmcsc10 scaled\magstep1
\font\xiiptbfs=cmb10  scaled\magstep1 \font\xiiptbmit=cmmib10
scaled\magstep1

\skewchar\xiiptmit='177 \skewchar\xiiptbmit='177 \skewchar\xiiptsy='60
\fontdimen16 \xiiptsy=\the\fontdimen17 \xiiptsy

\def\xiipt{\ifmmode\err@badsizechange\else
     \@mathfontinit
     \textfont0=\xiiptrm  \scriptfont0=\viiiptrm  \scriptscriptfont0=\viptrm
     \textfont1=\xiiptmit \scriptfont1=\viiiptmit \scriptscriptfont1=\viptmit
     \textfont2=\xiiptsy  \scriptfont2=\viiiptsy  \scriptscriptfont2=\viptsy
     \textfont3=\xiiptex  \scriptfont3=\xptex     \scriptscriptfont3=\xptex
     \textfont\itfam=\xiiptit
     \scriptfont\itfam=\viiiptit
     \scriptscriptfont\itfam=\viiptit
     \textfont\bffam=\xiiptbf
     \scriptfont\bffam=\viiiptbf
     \scriptscriptfont\bffam=\viptbf
     \textfont\bfsfam=\xiiptbfs
     \scriptfont\bfsfam=\viiiptbf
     \scriptscriptfont\bfsfam=\viptbf
     \textfont\bmitfam=\xiiptbmit
     \scriptfont\bmitfam=\viiiptmit
     \scriptscriptfont\bmitfam=\viptmit
     \@fontstyleinit
     \def\rm{\xiiptrm\fam=\z@}%
     \def\it{\xiiptit\fam=\itfam}%
     \def\sl{\xiiptsl}%
     \def\bf{\xiiptbf\fam=\bffam}%
     \def\tt{\xiipttt}%
     \def\ss{\xiiptss}%
     \def\sc{\xiiptsc}%
     \def\bfs{\xiiptbfs\fam=\bfsfam}%
     \def\bmit{\fam=\bmitfam}%
     \def\oldstyle{\xiiptmit\fam=\@ne}%
     \rm\fi}


\def\getxiiipt{%
     \font\xiiiptrm=cmr12  scaled\magstephalf
     \font\xiiiptmit=cmmi12 scaled\magstephalf
     \font\xiiiptsy=cmsy9  scaled\magstep2
     \font\xiiiptit=cmti12 scaled\magstephalf
     \font\xiiiptsl=cmsl12 scaled\magstephalf
     \font\xiiiptbf=cmbx12 scaled\magstephalf
     \font\xiiipttt=cmtt12 scaled\magstephalf
     \font\xiiiptss=cmss12 scaled\magstephalf
     \skewchar\xiiiptmit='177 \skewchar\xiiiptsy='60
     \fontdimen16 \xiiiptsy=\the\fontdimen17 \xiiiptsy}

\def\xiiipt{\ifmmode\err@badsizechange\else
     \@mathfontinit
     \textfont0=\xiiiptrm  \scriptfont0=\xptrm  \scriptscriptfont0=\viiptrm
     \textfont1=\xiiiptmit \scriptfont1=\xptmit \scriptscriptfont1=\viiptmit
     \textfont2=\xiiiptsy  \scriptfont2=\xptsy  \scriptscriptfont2=\viiptsy
     \textfont3=\xivptex   \scriptfont3=\xptex  \scriptscriptfont3=\xptex
     \textfont\itfam=\xiiiptit
     \scriptfont\itfam=\xptit
     \scriptscriptfont\itfam=\viiptit
     \textfont\bffam=\xiiiptbf
     \scriptfont\bffam=\xptbf
     \scriptscriptfont\bffam=\viiptbf
     \@fontstyleinit
     \def\rm{\xiiiptrm\fam=\z@}%
     \def\it{\xiiiptit\fam=\itfam}%
     \def\sl{\xiiiptsl}%
     \def\bf{\xiiiptbf\fam=\bffam}%
     \def\tt{\xiiipttt}%
     \def\ss{\xiiiptss}%
     \def\oldstyle{\xiiiptmit\fam=\@ne}%
     \rm\fi}


\font\xivptrm=cmr12   scaled\magstep1 \font\xivptmit=cmmi12
scaled\magstep1 \font\xivptsy=cmsy10  scaled\magstep2
\font\xivptex=cmex10  scaled\magstep2 \font\xivptit=cmti12
scaled\magstep1 \font\xivptsl=cmsl12  scaled\magstep1
\font\xivptbf=cmbx12  scaled\magstep1
\font\xivptss=cmss12  scaled\magstep1 \font\xivptsc=cmcsc10
scaled\magstep2 \font\xivptbfs=cmb10  scaled\magstep2
\font\xivptbmit=cmmib10 scaled\magstep2

\skewchar\xivptmit='177 \skewchar\xivptbmit='177 \skewchar\xivptsy='60
\fontdimen16 \xivptsy=\the\fontdimen17 \xivptsy

\def\xivpt{\ifmmode\err@badsizechange\else
     \@mathfontinit
     \textfont0=\xivptrm  \scriptfont0=\xptrm  \scriptscriptfont0=\viiptrm
     \textfont1=\xivptmit \scriptfont1=\xptmit \scriptscriptfont1=\viiptmit
     \textfont2=\xivptsy  \scriptfont2=\xptsy  \scriptscriptfont2=\viiptsy
     \textfont3=\xivptex  \scriptfont3=\xptex  \scriptscriptfont3=\xptex
     \textfont\itfam=\xivptit
     \scriptfont\itfam=\xptit
     \scriptscriptfont\itfam=\viiptit
     \textfont\bffam=\xivptbf
     \scriptfont\bffam=\xptbf
     \scriptscriptfont\bffam=\viiptbf
     \textfont\bfsfam=\xivptbfs
     \scriptfont\bfsfam=\xptbfs
     \scriptscriptfont\bfsfam=\viiptbf
     \textfont\bmitfam=\xivptbmit
     \scriptfont\bmitfam=\xptbmit
     \scriptscriptfont\bmitfam=\viiptmit
     \@fontstyleinit
     \def\rm{\xivptrm\fam=\z@}%
     \def\it{\xivptit\fam=\itfam}%
     \def\sl{\xivptsl}%
     \def\bf{\xivptbf\fam=\bffam}%
     \def\tt{\xivpttt}%
     \def\ss{\xivptss}%
     \def\sc{\xivptsc}%
     \def\bfs{\xivptbfs\fam=\bfsfam}%
     \def\bmit{\fam=\bmitfam}%
     \def\oldstyle{\xivptmit\fam=\@ne}%
     \rm\fi}


\font\xviiptrm=cmr17 \font\xviiptmit=cmmi12 scaled\magstep2
\font\xviiptsy=cmsy10 scaled\magstep3 \font\xviiptex=cmex10
scaled\magstep3 \font\xviiptit=cmti12 scaled\magstep2
\font\xviiptbf=cmbx12 scaled\magstep2 \font\xviiptbfs=cmb10
scaled\magstep3

\skewchar\xviiptmit='177 \skewchar\xviiptsy='60 \fontdimen16
\xviiptsy=\the\fontdimen17 \xviiptsy

\def\xviipt{\ifmmode\err@badsizechange\else
     \@mathfontinit
     \textfont0=\xviiptrm  \scriptfont0=\xiiptrm  \scriptscriptfont0=\viiiptrm
     \textfont1=\xviiptmit \scriptfont1=\xiiptmit \scriptscriptfont1=\viiiptmit
     \textfont2=\xviiptsy  \scriptfont2=\xiiptsy  \scriptscriptfont2=\viiiptsy
     \textfont3=\xviiptex  \scriptfont3=\xiiptex  \scriptscriptfont3=\xptex
     \textfont\itfam=\xviiptit
     \scriptfont\itfam=\xiiptit
     \scriptscriptfont\itfam=\viiiptit
     \textfont\bffam=\xviiptbf
     \scriptfont\bffam=\xiiptbf
     \scriptscriptfont\bffam=\viiiptbf
     \textfont\bfsfam=\xviiptbfs
     \scriptfont\bfsfam=\xiiptbfs
     \scriptscriptfont\bfsfam=\viiiptbf
     \@fontstyleinit
     \def\rm{\xviiptrm\fam=\z@}%
     \def\it{\xviiptit\fam=\itfam}%
     \def\bf{\xviiptbf\fam=\bffam}%
     \def\bfs{\xviiptbfs\fam=\bfsfam}%
     \def\oldstyle{\xviiptmit\fam=\@ne}%
     \rm\fi}


\font\xxiptrm=cmr17  scaled\magstep1


\def\xxipt{\ifmmode\err@badsizechange\else
     \@mathfontinit
     \@fontstyleinit
     \def\rm{\xxiptrm\fam=\z@}%
     \rm\fi}


\font\xxvptrm=cmr17  scaled\magstep2


\def\xxvpt{\ifmmode\err@badsizechange\else
     \@mathfontinit
     \@fontstyleinit
     \def\rm{\xxvptrm\fam=\z@}%
     \rm\fi}




\message{Loading jyTeX macros...}

\message{modifications to plain.tex,}


\def\newcount{\alloc@0\count\countdef\insc@unt}
\def\newdimen{\alloc@1\dimen\dimendef\insc@unt}
\def\newskip{\alloc@2\skip\skipdef\insc@unt}
\def\newmuskip{\alloc@3\muskip\muskipdef\@cclvi}
\def\newbox{\alloc@4\box\chardef\insc@unt}
\def\newtoks{\alloc@5\toks\toksdef\@cclvi}
\def\newhelp#1#2{\newtoks#1\global#1\expandafter{\csname#2\endcsname}}
\def\newread{\alloc@6\read\chardef\sixt@@n}
\def\newwrite{\alloc@7\write\chardef\sixt@@n}
\def\newfam{\alloc@8\fam\chardef\sixt@@n}
\def\newinsert#1{\global\advance\insc@unt by\m@ne
     \ch@ck0\insc@unt\count
     \ch@ck1\insc@unt\dimen
     \ch@ck2\insc@unt\skip
     \ch@ck4\insc@unt\box
     \allocationnumber=\insc@unt
     \global\chardef#1=\allocationnumber
     \wlog{\string#1=\string\insert\the\allocationnumber}}
\def\newif#1{\count@\escapechar \escapechar\m@ne
     \expandafter\expandafter\expandafter
          \xdef\@if#1{true}{\let\noexpand#1=\noexpand\iftrue}%
     \expandafter\expandafter\expandafter
          \xdef\@if#1{false}{\let\noexpand#1=\noexpand\iffalse}%
     \global\@if#1{false}\escapechar=\count@}


\newlinechar=`\^^J
\overfullrule=0pt




\let\itfam=\undefined

\let\bffam=\undefined

\count18=3


\chardef\sharps="19


\mathchardef\alpha="710B \mathchardef\beta="710C \mathchardef\gamma="710D
\mathchardef\delta="710E \mathchardef\epsilon="710F
\mathchardef\zeta="7110 \mathchardef\eta="7111 \mathchardef\theta="7112
\mathchardef\iota="7113 \mathchardef\kappa="7114
\mathchardef\lambda="7115 \mathchardef\mu="7116 \mathchardef\nu="7117
\mathchardef\xi="7118 \mathchardef\pi="7119 \mathchardef\rho="711A
\mathchardef\sigma="711B \mathchardef\tau="711C
\mathchardef\upsilon="711D \mathchardef\phi="711E \mathchardef\chi="711F
\mathchardef\psi="7120 \mathchardef\omega="7121
\mathchardef\varepsilon="7122 \mathchardef\vartheta="7123
\mathchardef\varpi="7124 \mathchardef\varrho="7125
\mathchardef\varsigma="7126 \mathchardef\varphi="7127
\mathchardef\imath="717B \mathchardef\jmath="717C \mathchardef\ell="7160
\mathchardef\wp="717D \mathchardef\partial="7140 \mathchardef\flat="715B
\mathchardef\natural="715C \mathchardef\sharp="715D



\def\angle{{\vbox{\ialign{$\m@th\scriptstyle##$\crcr
     \not\mathrel{\mkern14mu}\crcr
     \noalign{\nointerlineskip}
     \mkern2.5mu\leaders\hrule height.34\rp@\hfill\mkern2.5mu\crcr}}}}
\def\vdots{\vbox{\baselineskip4\rp@ \lineskiplimit\z@
     \kern6\rp@\hbox{.}\hbox{.}\hbox{.}}}
\def\ddots{\mathinner{\mkern1mu\raise7\rp@\vbox{\kern7\rp@\hbox{.}}\mkern2mu
     \raise4\rp@\hbox{.}\mkern2mu\raise\rp@\hbox{.}\mkern1mu}}
\def\overrightarrow#1{\vbox{\ialign{##\crcr
     \rightarrowfill\crcr
     \noalign{\kern-\rp@\nointerlineskip}
     $\hfil\displaystyle{#1}\hfil$\crcr}}}
\def\overleftarrow#1{\vbox{\ialign{##\crcr
     \leftarrowfill\crcr
     \noalign{\kern-\rp@\nointerlineskip}
     $\hfil\displaystyle{#1}\hfil$\crcr}}}
\def\overbrace#1{\mathop{\vbox{\ialign{##\crcr
     \noalign{\kern3\rp@}
     \downbracefill\crcr
     \noalign{\kern3\rp@\nointerlineskip}
     $\hfil\displaystyle{#1}\hfil$\crcr}}}\limits}
\def\underbrace#1{\mathop{\vtop{\ialign{##\crcr
     $\hfil\displaystyle{#1}\hfil$\crcr
     \noalign{\kern3\rp@\nointerlineskip}
     \upbracefill\crcr
     \noalign{\kern3\rp@}}}}\limits}
\def\big#1{{\hbox{$\left#1\vbox to8.5\rp@ {}\right.\n@space$}}}
\def\Big#1{{\hbox{$\left#1\vbox to11.5\rp@ {}\right.\n@space$}}}
\def\bigg#1{{\hbox{$\left#1\vbox to14.5\rp@ {}\right.\n@space$}}}
\def\Bigg#1{{\hbox{$\left#1\vbox to17.5\rp@ {}\right.\n@space$}}}
\def\@vereq#1#2{\lower.5\rp@\vbox{\baselineskip\z@skip\lineskip-.5\rp@
     \ialign{$\m@th#1\hfil##\hfil$\crcr#2\crcr=\crcr}}}
\def\rlh@#1{\vcenter{\hbox{\ooalign{\raise2\rp@
     \hbox{$#1\rightharpoonup$}\crcr
     $#1\leftharpoondown$}}}}
\def\bordermatrix#1{\begingroup\m@th
     \setbox\z@\vbox{%
          \def\cr{\crcr\noalign{\kern2\rp@\global\let\cr\endline}}%
          \ialign{$##$\hfil\kern2\rp@\kern\p@renwd
               &\thinspace\hfil$##$\hfil&&\quad\hfil$##$\hfil\crcr
               \omit\strut\hfil\crcr
               \noalign{\kern-\baselineskip}%
               #1\crcr\omit\strut\cr}}%
     \setbox\tw@\vbox{\unvcopy\z@\global\setbox\@ne\lastbox}%
     \setbox\tw@\hbox{\unhbox\@ne\unskip\global\setbox\@ne\lastbox}%
     \setbox\tw@\hbox{$\kern\wd\@ne\kern-\p@renwd\left(\kern-\wd\@ne
          \global\setbox\@ne\vbox{\box\@ne\kern2\rp@}%
          \vcenter{\kern-\ht\@ne\unvbox\z@\kern-\baselineskip}%
          \,\right)$}%
     \null\;\vbox{\kern\ht\@ne\box\tw@}\endgroup}
\def\endinsert{\egroup
     \if@mid\dimen@\ht\z@
          \advance\dimen@\dp\z@
          \advance\dimen@12\rp@
          \advance\dimen@\pagetotal
          \ifdim\dimen@>\pagegoal\@midfalse\p@gefalse\fi
     \fi
     \if@mid\bigskip\box\z@
          \bigbreak
     \else\insert\topins{\penalty100 \splittopskip\z@skip
               \splitmaxdepth\maxdimen\floatingpenalty\z@
               \ifp@ge\dimen@\dp\z@
                    \vbox to\vsize{\unvbox\z@\kern-\dimen@}%
               \else\box\z@\nobreak\bigskip
               \fi}%
     \fi
     \endgroup}


\def\cases#1{\left\{\,\vcenter{\m@th
     \ialign{$##\hfil$&\quad##\hfil\crcr#1\crcr}}\right.}
\def\matrix#1{\null\,\vcenter{\m@th
     \ialign{\hfil$##$\hfil&&\quad\hfil$##$\hfil\crcr
          \mathstrut\crcr
          \noalign{\kern-\baselineskip}
          #1\crcr
          \mathstrut\crcr
          \noalign{\kern-\baselineskip}}}\,}


\newif\ifraggedbottom

\def\raggedbottom{\ifraggedbottom\else
     \advance\topskip by\z@ plus60pt \raggedbottomtrue\fi}%
\def\normalbottom{\ifraggedbottom
     \advance\topskip by\z@ plus-60pt \raggedbottomfalse\fi}

\message{hacks,}


\toksdef\toks@i=1 \toksdef\toks@ii=2


\def\TeX{T\kern-.1667em \lower.5ex \hbox{E}\kern-.125em X\null}
\def\jyTeX{{\leavevmode
     \raise.587ex \hbox{\it\j}\kern-.1em \lower.048ex \hbox{\it y}\kern-.12em
     \TeX}}

\let\then=\iftrue
\def\ifnoarg#1\then{\def\hack@{#1}\ifx\hack@\empty}
\def\ifundefined#1\then{%
     \expandafter\ifx\csname\expandafter\blank\string#1\endcsname\relax}
\def\useif#1\then{\csname#1\endcsname}
\def\usename#1{\csname#1\endcsname}
\def\useafter#1#2{\expandafter#1\csname#2\endcsname}

\long\def\loop#1\repeat{\def\@iterate{#1\expandafter\@iterate\fi}\@iterate
     \let\@iterate=\relax}

\let\TeXend=\end
\def\begin#1{\begingroup\def\@@blockname{#1}\usename{begin#1}}
\def\end#1{\usename{end#1}\def\hack@{#1}%
     \ifx\@@blockname\hack@
          \endgroup
     \else\err@badgroup\hack@\@@blockname
     \fi}
\def\@@blockname{}

\def\defaultoption[#1]#2{%
     \def\hack@{\ifx\hack@ii[\toks@={#2}\else\toks@={#2[#1]}\fi\the\toks@}%
     \futurelet\hack@ii\hack@}

\def\markup#1{\let\@@marksf=\empty
     \ifhmode\edef\@@marksf{\spacefactor=\the\spacefactor\relax}\/\fi
     ${}^{\hbox{\subscriptfonts#1}}$\@@marksf}


\newtoks\shortyear
\newtoks\militaryhour
\newtoks\standardhour
\newtoks\minute
\newtoks\amorpm

\def\settime{\count@=\time\divide\count@ by60
     \militaryhour=\expandafter{\number\count@}%
     {\multiply\count@ by-60 \advance\count@ by\time
          \xdef\hack@{\ifnum\count@<10 0\fi\number\count@}}%
     \minute=\expandafter{\hack@}%
     \ifnum\count@<12
          \amorpm={am}
     \else\amorpm={pm}
          \ifnum\count@>12 \advance\count@ by-12 \fi
     \fi
     \standardhour=\expandafter{\number\count@}%
     \def\hack@19##1##2{\shortyear={##1##2}}%
          \expandafter\hack@\the\year}

\def\monthword#1{%
     \ifcase#1
          $\bullet$\err@badcountervalue{monthword}%
          \or January\or February\or March\or April\or May\or June%
          \or July\or August\or September\or October\or November\or December%
     \else$\bullet$\err@badcountervalue{monthword}%
     \fi}

\def\monthabbr#1{%
     \ifcase#1
          $\bullet$\err@badcountervalue{monthabbr}%
          \or Jan\or Feb\or Mar\or Apr\or May\or Jun%
          \or Jul\or Aug\or Sep\or Oct\or Nov\or Dec%
     \else$\bullet$\err@badcountervalue{monthabbr}%
     \fi}

\def\militarytime{\the\militaryhour:\the\minute}
\def\standardtime{\the\standardhour:\the\minute}


\def\@setnumstyle#1#2{\expandafter\global\expandafter\expandafter
     \expandafter\let\expandafter\expandafter
     \csname @\expandafter\blank\string#1style\endcsname
     \csname#2\endcsname}
\def\numstyle#1{\usename{@\expandafter\blank\string#1style}#1}
\def\ifblank#1\then{\useafter\ifx{@\expandafter\blank\string#1}\blank}

\def\blank#1{}

\def\Roman#1{\expandafter\uppercase\expandafter{\romannumeral#1}}
\def\alphabetic#1{%
     \ifcase#1
          $\bullet$\err@badcountervalue{alphabetic}%
          \or a\or b\or c\or d\or e\or f\or g\or h\or i\or j\or k\or l\or m%
          \or n\or o\or p\or q\or r\or s\or t\or u\or v\or w\or x\or y\or z%
     \else$\bullet$\err@badcountervalue{alphabetic}%
     \fi}
\def\Alphabetic#1{\expandafter\uppercase\expandafter{\alphabetic{#1}}}
\def\symbols#1{%
     \ifcase#1
          $\bullet$\err@badcountervalue{symbols}%
          \or*\or\dag\or\ddag\or\S\or$\|$%
          \or**\or\dag\dag\or\ddag\ddag\or\S\S\or$\|\|$%
     \else$\bullet$\err@badcountervalue{symbols}%
     \fi}


\catcode`\^^?=13 \def^^?{\relax}

\def\trimleading#1\to#2{\edef#2{#1}%
     \expandafter\@trimleading\expandafter#2#2^^?^^?}
\def\@trimleading#1#2#3^^?{\ifx#2^^?\def#1{}\else\def#1{#2#3}\fi}

\def\trimtrailing#1\to#2{\edef#2{#1}%
     \expandafter\@trimtrailing\expandafter#2#2^^? ^^?\relax}
\def\@trimtrailing#1#2 ^^?#3{\ifx#3\relax\toks@={}%
     \else\def#1{#2}\toks@={\trimtrailing#1\to#1}\fi
     \the\toks@}

\def\trim#1\to#2{\trimleading#1\to#2\trimtrailing#2\to#2}

\catcode`\^^?=15


\long\def\additemL#1\to#2{\toks@={\^^\{#1}}\toks@ii=\expandafter{#2}%
     \xdef#2{\the\toks@\the\toks@ii}}

\long\def\additemR#1\to#2{\toks@={\^^\{#1}}\toks@ii=\expandafter{#2}%
     \xdef#2{\the\toks@ii\the\toks@}}

\def\getitemL#1\to#2{\expandafter\@getitemL#1\hack@#1#2}
\def\@getitemL\^^\#1#2\hack@#3#4{\def#4{#1}\def#3{#2}}

\message{font macros,}


\newdimen\rp@
\newcount\@@sizeindex \@@sizeindex=0
\newcount\@@factori
\newcount\@@factorii
\newcount\@@factoriii
\newcount\@@factoriv

\countdef\maxfam=18
\newfam\itfam
\newfam\bffam
\newfam\bfsfam
\newfam\bmitfam

\def\@mathfontinit{\count@=4
     \loop\textfont\count@=\nullfont
          \scriptfont\count@=\nullfont
          \scriptscriptfont\count@=\nullfont
          \ifnum\count@<\maxfam\advance\count@ by\@ne
     \repeat}

\def\@fontstyleinit{%
     \def\it{\err@fontnotavailable\it}%
     \def\bf{\err@fontnotavailable\bf}%
     \def\bfs{\err@bfstobf}%
     \def\bmit{\err@fontnotavailable\bmit}%
     \def\sc{\err@fontnotavailable\sc}%
     \def\sl{\err@sltoit}%
     \def\ss{\err@fontnotavailable\ss}%
     \def\tt{\err@fontnotavailable\tt}}

\def\@parameterinit#1{\rm\rp@=.1em \@getscaling{#1}%
     \let\^^\=\@doscaling\scalingskipslist
     \setbox\strutbox=\hbox{\vrule
          height.708\baselineskip depth.292\baselineskip width\z@}}

\def\@getfactor#1#2#3#4{\@@factori=#1 \@@factorii=#2
     \@@factoriii=#3 \@@factoriv=#4}

\def\@getscaling#1{\count@=#1 \advance\count@ by-\@@sizeindex\@@sizeindex=#1
     \ifnum\count@<0
          \let\@mulordiv=\divide
          \let\@divormul=\multiply
          \multiply\count@ by\m@ne
     \else\let\@mulordiv=\multiply
          \let\@divormul=\divide
     \fi
     \edef\@@scratcha{\ifcase\count@                {1}{1}{1}{1}\or
          {1}{7}{23}{3}\or     {2}{5}{3}{1}\or      {9}{89}{13}{1}\or
          {6}{25}{6}{1}\or     {8}{71}{14}{1}\or    {6}{25}{36}{5}\or
          {1}{7}{53}{4}\or     {12}{125}{108}{5}\or {3}{14}{53}{5}\or
          {6}{41}{17}{1}\or    {13}{31}{13}{2}\or   {9}{107}{71}{2}\or
          {11}{139}{124}{3}\or {1}{6}{43}{2}\or     {10}{107}{42}{1}\or
          {1}{5}{43}{2}\or     {5}{69}{65}{1}\or    {11}{97}{91}{2}\fi}%
     \expandafter\@getfactor\@@scratcha}

\def\@doscaling#1{\@mulordiv#1by\@@factori\@divormul#1by\@@factorii
     \@mulordiv#1by\@@factoriii\@divormul#1by\@@factoriv}


\newskip\headskip
\newskip\footskip

\def\typesize=#1pt{\count@=#1 \advance\count@ by-10
     \ifcase\count@
          \@setsizex\or\err@badtypesize\or
          \@setsizexii\or\err@badtypesize\or
          \@setsizexiv
     \else\err@badtypesize
     \fi}

\def\@setsizex{\getixpt
     \def\subsubscriptfonts{\vpt}%
          \def\subsubscriptsize{\vpt\@parameterinit{-8}}%
     \def\subscriptfonts{\viipt}\def\subscriptsize{\viipt\@parameterinit{-4}}%
     \def\footnotefonts{\viiipt}\def\footnotesize{\viiipt\@parameterinit{-2}}%
     \def\smallfonts{\ixpt}\def\smallsize{\ixpt\@parameterinit{-1}}%
     \def\normalfonts{\xpt}\def\normalsize{\xpt\@parameterinit{0}}%
     \def\bigfonts{\xiipt}\def\bigsize{\xiipt\@parameterinit{2}}%
     \def\Bigfonts{\xivpt}\def\Bigsize{\xivpt\@parameterinit{4}}%
     \def\biggfonts{\xviipt}\def\biggsize{\xviipt\@parameterinit{6}}%
     \def\Biggfonts{\xxipt}\def\Biggsize{\xxipt\@parameterinit{8}}%
     \def\tinyfonts{\vpt}\def\tinysize{\vpt\@parameterinit{-8}}%
     \def\HUGEFONTS{\xxvpt}\def\HUGESIZE{\xxvpt\@parameterinit{10}}%
     \normalsize\fixedskipslist}

\def\@setsizexii{\getxipt
     \def\subsubscriptfonts{\vipt}%
          \def\subsubscriptsize{\vipt\@parameterinit{-6}}%
     \def\subscriptfonts{\viiipt}%
          \def\subscriptsize{\viiipt\@parameterinit{-2}}%
     \def\footnotefonts{\xpt}\def\footnotesize{\xpt\@parameterinit{0}}%
     \def\smallfonts{\xipt}\def\smallsize{\xipt\@parameterinit{1}}%
     \def\normalfonts{\xiipt}\def\normalsize{\xiipt\@parameterinit{2}}%
     \def\bigfonts{\xivpt}\def\bigsize{\xivpt\@parameterinit{4}}%
     \def\Bigfonts{\xviipt}\def\Bigsize{\xviipt\@parameterinit{6}}%
     \def\biggfonts{\xxipt}\def\biggsize{\xxipt\@parameterinit{8}}%
     \def\Biggfonts{\xxvpt}\def\Biggsize{\xxvpt\@parameterinit{10}}%
     \def\tinyfonts{\vpt}\def\tinysize{\vpt\@parameterinit{-8}}%
     \def\HUGEFONTS{\xxvpt}\def\HUGESIZE{\xxvpt\@parameterinit{10}}%
     \normalsize\fixedskipslist}

\def\@setsizexiv{\getxiiipt
     \def\subsubscriptfonts{\viipt}%
          \def\subsubscriptsize{\viipt\@parameterinit{-4}}%
     \def\subscriptfonts{\xpt}\def\subscriptsize{\xpt\@parameterinit{0}}%
     \def\footnotefonts{\xiipt}\def\footnotesize{\xiipt\@parameterinit{2}}%
     \def\smallfonts{\xiiipt}\def\smallsize{\xiiipt\@parameterinit{3}}%
     \def\normalfonts{\xivpt}\def\normalsize{\xivpt\@parameterinit{4}}%
     \def\bigfonts{\xviipt}\def\bigsize{\xviipt\@parameterinit{6}}%
     \def\Bigfonts{\xxipt}\def\Bigsize{\xxipt\@parameterinit{8}}%
     \def\biggfonts{\xxvpt}\def\biggsize{\xxvpt\@parameterinit{10}}%
     \def\Biggfonts{\err@sizetoolarge\Biggfonts\HUGEFONTS}%
          \def\Biggsize{\err@sizetoolarge\Biggsize\HUGESIZE}%
     \def\tinyfonts{\vpt}\def\tinysize{\vpt\@parameterinit{-8}}%
     \def\HUGEFONTS{\xxvpt}\def\HUGESIZE{\xxvpt\@parameterinit{10}}%
     \normalsize\fixedskipslist}

\def\subsubscriptfonts{\vpt} \def\subsubscriptsize{\vpt\@parameterinit{-8}}
\def\subscriptfonts{\viipt}  \def\subscriptsize{\viipt\@parameterinit{-4}}
\def\footnotefonts{\viiipt}  \def\footnotesize{\viiipt\@parameterinit{-2}}
\def\smallfonts{\err@sizenotavailable\smallfonts}
                             \def\smallsize{\ixpt\@parameterinit{-1}}
\def\normalfonts{\xpt}       \def\normalsize{\xpt\@parameterinit{0}}
\def\bigfonts{\xiipt}        \def\bigsize{\xiipt\@parameterinit{2}}
\def\Bigfonts{\xivpt}        \def\Bigsize{\xivpt\@parameterinit{4}}
\def\biggfonts{\xviipt}      \def\biggsize{\xviipt\@parameterinit{6}}
\def\Biggfonts{\xxipt}       \def\Biggsize{\xxipt\@parameterinit{8}}
\def\tinyfonts{\vpt}         \def\tinysize{\vpt\@parameterinit{-8}}
\def\HUGEFONTS{\xxvpt}       \def\HUGESIZE{\xxvpt\@parameterinit{10}}

\message{document layout,}


\newtoks\everyoutput \everyoutput={}
\newdimen\depthofpage
\newcount\pagenum \pagenum=0

\newdimen\oddtopmargin  \newdimen\eventopmargin
\newdimen\oddleftmargin \newdimen\evenleftmargin
\newtoks\oddhead        \newtoks\evenhead
\newtoks\oddfoot        \newtoks\evenfoot

\def\topmargin{\afterassignment\@seteventop\oddtopmargin}
\def\leftmargin{\afterassignment\@setevenleft\oddleftmargin}
\def\head{\afterassignment\@setevenhead\oddhead}
\def\foot{\afterassignment\@setevenfoot\oddfoot}

\def\@seteventop{\eventopmargin=\oddtopmargin}
\def\@setevenleft{\evenleftmargin=\oddleftmargin}
\def\@setevenhead{\evenhead=\oddhead}
\def\@setevenfoot{\evenfoot=\oddfoot}

\def\pagenumstyle#1{\@setnumstyle\pagenum{#1}}

\newif\ifdraft
\def\draft{\drafttrue\leftmargin=.5in \overfullrule=5pt }

\def\outputstyle#1{\global\expandafter\let\expandafter
          \@outputstyle\csname#1output\endcsname
     \usename{#1setup}}

\output={\@outputstyle}

\def\normaloutput{\the\everyoutput
     \global\advance\pagenum by\@ne
     \ifodd\pagenum
          \voffset=\oddtopmargin \hoffset=\oddleftmargin
     \else\voffset=\eventopmargin \hoffset=\evenleftmargin
     \fi
     \advance\voffset by-1in  \advance\hoffset by-1in
     \count0=\pagenum
     \expandafter\shipout\pagebox
     \ifnum\outputpenalty>-\@MM\else\dosupereject\fi}

\newdimen\fullhsize
\newbox\leftpage
\newcount\leftpagenum
\newcount\outputpagenum \outputpagenum=0
\let\leftorright=L

\def\twoupoutput{\the\everyoutput
     \global\advance\pagenum by\@ne
     \if L\leftorright
          \global\setbox\leftpage=\leftline{\pagebox}%
          \global\leftpagenum=\pagenum
          \global\let\leftorright=R%
     \else\global\advance\outputpagenum by\@ne
          \ifodd\outputpagenum
               \voffset=\oddtopmargin \hoffset=\oddleftmargin
          \else\voffset=\eventopmargin \hoffset=\evenleftmargin
          \fi
          \advance\voffset by-1in  \advance\hoffset by-1in
          \count0=\leftpagenum \count1=\pagenum
          \shipout\vbox{\hbox to\fullhsize
               {\box\leftpage\hfil\leftline{\pagebox}}}%
          \global\let\leftorright=L%
     \fi
     \ifnum\outputpenalty>-\@MM
     \else\dosupereject
          \if R\leftorright
               \globaldefs=\@ne\head={\hfil}\foot={\hfil}\globaldefs=\z@
               \null\newpage
          \fi
     \fi}

\def\pagebox{\vbox{\makeheadline\pagebody\makefootline}}

\def\makeheadline{%
     \vbox to\z@{\baselinestretch=\@m
          \vskip\topskip\vskip-.708\baselineskip\vskip-\headskip
          \line{\vbox to\ht\strutbox{}%
               \ifodd\pagenum\the\oddhead\else\the\evenhead\fi}%
          \vss}%
     \nointerlineskip}

\def\pagebody{\vbox to\vsize{%
     \boxmaxdepth\maxdepth
     \ifvoid\topins\else\unvbox\topins\fi
     \depthofpage=\dp255
     \unvbox255
     \ifraggedbottom\kern-\depthofpage\vfil\fi
     \ifvoid\footins
     \else\vskip\skip\footins
          \footnoterule
          \unvbox\footins
          \vskip-\footnoteskip
     \fi}}

\def\makefootline{\baselineskip=\footskip
     \line{\ifodd\pagenum\the\oddfoot\else\the\evenfoot\fi}}


\newskip\abovechapterskip
\newskip\belowchapterskip
\newskip\abovesectionskip
\newskip\belowsectionskip
\newskip\abovesubsectionskip
\newskip\belowsubsectionskip

\def\chapterstyle#1{\global\expandafter\let\expandafter\@chapterstyle
     \csname#1text\endcsname}
\def\sectionstyle#1{\global\expandafter\let\expandafter\@sectionstyle
     \csname#1text\endcsname}
\def\subsectionstyle#1{\global\expandafter\let\expandafter\@subsectionstyle
     \csname#1text\endcsname}

\def\chapter#1{%
     \ifdim\lastskip=17sp \else\chapterbreak\vskip\abovechapterskip\fi
     \@chapterstyle{\ifblank\chapternumstyle\then
          \else\newchapternum=\next\chapternumformat\ \fi#1}%
     \nobreak\vskip\belowchapterskip\vskip17sp }

\def\section#1{%
     \ifdim\lastskip=17sp \else\sectionbreak\vskip\abovesectionskip\fi
     \@sectionstyle{\ifblank\sectionnumstyle\then
          \else\newsectionnum=\next\sectionnumformat\ \fi#1}%
     \nobreak\vskip\belowsectionskip\vskip17sp }

\def\subsection#1{%
     \ifdim\lastskip=17sp \else\subsectionbreak\vskip\abovesubsectionskip\fi
     \@subsectionstyle{\ifblank\subsectionnumstyle\then
          \else\newsubsectionnum=\next\subsectionnumformat\ \fi#1}%
     \nobreak\vskip\belowsubsectionskip\vskip17sp }


\let\TeXunderline=\underline
\let\TeXoverline=\overline
\def\underline#1{\relax\ifmmode\TeXunderline{#1}\else
     $\TeXunderline{\hbox{#1}}$\fi}
\def\overline#1{\relax\ifmmode\TeXoverline{#1}\else
     $\TeXoverline{\hbox{#1}}$\fi}

\def\baselinestretch{\afterassignment\@baselinestretch\count@}
\def\@baselinestretch{\baselineskip=\normalbaselineskip
     \divide\baselineskip by\@m\baselineskip=\count@\baselineskip
     \setbox\strutbox=\hbox{\vrule
          height.708\baselineskip depth.292\baselineskip width\z@}%
     \bigskipamount=\the\baselineskip
          plus.25\baselineskip minus.25\baselineskip
     \medskipamount=.5\baselineskip
          plus.125\baselineskip minus.125\baselineskip
     \smallskipamount=.25\baselineskip
          plus.0625\baselineskip minus.0625\baselineskip}

\def\\{\ifhmode\ifnum\lastpenalty=-\@M\else\hfil\penalty-\@M\fi\fi
     \ignorespaces}
\def\newpage{\vfil\break}

\def\lefttext#1{\par{\@text\leftskip=\z@\rightskip=\centering
     \noindent#1\par}}
\def\righttext#1{\par{\@text\leftskip=\centering\rightskip=\z@
     \noindent#1\par}}
\def\centertext#1{\par{\@text\leftskip=\centering\rightskip=\centering
     \noindent#1\par}}
\def\@text{\parindent=\z@ \parfillskip=\z@ \everypar={}%
     \spaceskip=.3333em \xspaceskip=.5em
     \def\\{\ifhmode\ifnum\lastpenalty=-\@M\else\penalty-\@M\fi\fi
          \ignorespaces}}

\def\beginleft{\par\@text\leftskip=\z@ \rightskip=\centering}
     
\def\beginright{\par\@text\leftskip=\centering\rightskip=\z@ }
     
\def\begincenter{\par\@text\leftskip=\centering\rightskip=\centering}

\def\beginnarrow{\defaultoption[\parindent]\@beginnarrow}
\def\@beginnarrow[#1]{\par\advance\leftskip by#1\advance\rightskip by#1}

\begingroup
\catcode`\[=1 \catcode`\{=11 \gdef\beginignore[\endgroup\bgroup
     \catcode`\e=0 \catcode`\\=12 \catcode`\{=11 \catcode`\f=12 \let\or=\relax
     \let\nd{ignor=\fi \let\}=\egroup
     \iffalse}
\endgroup

\long\def\marginnote#1{\leavevmode
     \edef\@marginsf{\spacefactor=\the\spacefactor\relax}%
     \ifdraft\strut\vadjust{%
          \hbox to\z@{\hskip\hsize\hskip.1in
               \vbox to\z@{\vskip-\dp\strutbox
                    \marginnoteformat
                    \vskip-\ht\strutbox
                    \noindent\strut#1\par
                    \vss}%
               \hss}}%
     \fi
     \@marginsf}


\newtoks\everybye \everybye={\par\vfil}
\outer\def\bye{\the\everybye
     \footnotecheck
     \prelabelcheck
     \streamcheck
     \supereject
     \TeXend}

\message{footnotes,}

\newcount\footnotenum \footnotenum=0
\newskip\footnoteskip
\let\@footnotelist=\empty

\def\footnotenumstyle#1{\@setnumstyle\footnotenum{#1}%
     \useafter\ifx{@footnotenumstyle}\symbols
          \global\let\@footup=\empty
     \else\global\let\@footup=\markup
     \fi}

\def\footnote{\footnotecheck\defaultoption[]\@footnote}
\def\@footnote[#1]{\@footnotemark[#1]\@footnotetext}

\def\footnotemark{\defaultoption[]\@footnotemark}
\def\@footnotemark[#1]{\let\@footsf=\empty
     \ifhmode\edef\@footsf{\spacefactor=\the\spacefactor\relax}\/\fi
     \ifnoarg#1\then
          \global\advance\footnotenum by\@ne
          \@footup{\footnotenumformat}%
          \edef\@@foota{\footnotenum=\the\footnotenum\relax}%
          \expandafter\additemR\expandafter\@footup\expandafter
               {\@@foota\footnotenumformat}\to\@footnotelist
          \global\let\@footnotelist=\@footnotelist
     \else\markup{#1}%
          \additemR\markup{#1}\to\@footnotelist
          \global\let\@footnotelist=\@footnotelist
     \fi
     \@footsf}

\def\footnotetext{%
     \ifx\@footnotelist\empty\err@extrafootnotetext\else\@footnotetext\fi}
\def\@footnotetext{%
     \getitemL\@footnotelist\to\@@foota
     \global\let\@footnotelist=\@footnotelist
     \insert\footins\bgroup
     \footnoteformat
     \splittopskip=\ht\strutbox\splitmaxdepth=\dp\strutbox
     \interlinepenalty=\interfootnotelinepenalty\floatingpenalty=\@MM
     \noindent\llap{\@@foota}\strut
     \bgroup\aftergroup\@footnoteend
     \let\@@scratcha=}
\def\@footnoteend{\strut\par\vskip\footnoteskip\egroup}

\def\footnoterule{\normalfonts
     \kern-.3em \hrule width2in height.04em \kern .26em }

\def\footnotecheck{%
     \ifx\@footnotelist\empty
     \else\err@extrafootnotemark
          \global\let\@footnotelist=\empty
     \fi}

\message{labels,}

\let\@@labeldef=\xdef
\newif\if@labelfile
\newwrite\@labelfile
\let\@prelabellist=\empty

\def\label#1#2{\trim#1\to\@@labarg\edef\@@labtext{#2}%
     \edef\@@labname{lab@\@@labarg}%
     \useafter\ifundefined\@@labname\then\else\@yeslab\fi
     \useafter\@@labeldef\@@labname{#2}%
     \ifstreaming
          \expandafter\toks@\expandafter\expandafter\expandafter
               {\csname\@@labname\endcsname}%
          \immediate\write\streamout{\noexpand\label{\@@labarg}{\the\toks@}}%
     \fi}
\def\@yeslab{%
     \useafter\ifundefined{if\@@labname}\then
          \err@labelredef\@@labarg
     \else\useif{if\@@labname}\then
               \err@labelredef\@@labarg
          \else\global\usename{\@@labname true}%
               \useafter\ifundefined{pre\@@labname}\then
               \else\useafter\ifx{pre\@@labname}\@@labtext
                    \else\err@badlabelmatch\@@labarg
                    \fi
               \fi
               \if@labelfile
               \else\global\@labelfiletrue
                    \immediate\write\sixt@@n{--> Creating file \jobname.lab}%
                    \immediate\openout\@labelfile=\jobname.lab
               \fi
               \immediate\write\@labelfile
                    {\noexpand\prelabel{\@@labarg}{\@@labtext}}%
          \fi
     \fi}

\def\putlab#1{\trim#1\to\@@labarg\edef\@@labname{lab@\@@labarg}%
     \useafter\ifundefined\@@labname\then\@nolab\else\usename\@@labname\fi}
\def\@nolab{%
     \useafter\ifundefined{pre\@@labname}\then
          \undefinedlabelformat
          \err@needlabel\@@labarg
          \useafter\xdef\@@labname{\undefinedlabelformat}%
     \else\usename{pre\@@labname}%
          \useafter\xdef\@@labname{\usename{pre\@@labname}}%
     \fi
     \useafter\newif{if\@@labname}%
     \expandafter\additemR\@@labarg\to\@prelabellist}

\def\prelabel#1{\useafter\gdef{prelab@#1}}

\def\ifundefinedlabel#1\then{%
     \expandafter\ifx\csname lab@#1\endcsname\relax}
\def\useiflab#1\then{\csname iflab@#1\endcsname}

\def\prelabelcheck{{%
     \def\^^\##1{\useiflab{##1}\then\else\err@undefinedlabel{##1}\fi}%
     \@prelabellist}}

\message{equation numbering,}

\newcount\chapternum
\newcount\sectionnum
\newcount\subsectionnum
\newcount\equationnum
\newcount\subequationnum
\newcount\figurenum
\newcount\subfigurenum
\newcount\tablenum
\newcount\subtablenum

\newif\if@subeqncount
\newif\if@subfigcount
\newif\if@subtblcount

\def\newchapternum{\newsectionnum=\z@\@resetnum\chapternum}
\def\newsectionnum{\newsubsectionnum=\z@\@resetnum\sectionnum}
\def\newsubsectionnum{\newequationnum=\z@\newfigurenum=\z@\newtablenum=\z@
     \@resetnum\subsectionnum}
\def\newequationnum{\newsubequationnum=\z@\@resetnum\equationnum}
\def\newsubequationnum{\@resetnum\subequationnum}
\def\newfigurenum{\newsubfigurenum=\z@\@resetnum\figurenum}
\def\newsubfigurenum{\@resetnum\subfigurenum}
\def\newtablenum{\newsubtablenum=\z@\@resetnum\tablenum}
\def\newsubtablenum{\@resetnum\subtablenum}

\def\@resetnum#1{\global\advance#1by1 \edef\next{\the#1\relax}\global#1}

\newchapternum=0

\def\chapternumstyle#1{\@setnumstyle\chapternum{#1}}
\def\sectionnumstyle#1{\@setnumstyle\sectionnum{#1}}
\def\subsectionnumstyle#1{\@setnumstyle\subsectionnum{#1}}
\def\equationnumstyle#1{\@setnumstyle\equationnum{#1}}
\def\subequationnumstyle#1{\@setnumstyle\subequationnum{#1}%
     \ifblank\subequationnumstyle\then\global\@subeqncountfalse\fi
     \ignorespaces}
\def\figurenumstyle#1{\@setnumstyle\figurenum{#1}}
\def\subfigurenumstyle#1{\@setnumstyle\subfigurenum{#1}%
     \ifblank\subfigurenumstyle\then\global\@subfigcountfalse\fi
     \ignorespaces}
\def\tablenumstyle#1{\@setnumstyle\tablenum{#1}}
\def\subtablenumstyle#1{\@setnumstyle\subtablenum{#1}%
     \ifblank\subtablenumstyle\then\global\@subtblcountfalse\fi
     \ignorespaces}

\def\eqnlabel#1{%
     \if@subeqncount
          \newsubequationnum=\next
     \else\newequationnum=\next
          \ifblank\subequationnumstyle\then
          \else\global\@subeqncounttrue
               \newsubequationnum=\@ne
          \fi
     \fi
     \label{#1}{\puteqnformat}(\puteqn{#1})%
     \ifdraft\rlap{\hskip.1in{\tt#1}}\fi}

\let\puteqn=\putlab

\def\equation#1#2{\useafter\gdef{eqn@#1}{#2\eqno\eqnlabel{#1}}}
\def\Equation#1{\useafter\gdef{eqn@#1}}

\def\putequation#1{\useafter\ifundefined{eqn@#1}\then
     \err@undefinedeqn{#1}\else\usename{eqn@#1}\fi}

\def\eqnseriesstyle#1{\gdef\@eqnseriesstyle{#1}}
\def\begineqnseries{\subequationnumstyle{\@eqnseriesstyle}%
     \defaultoption[]\@begineqnseries}
\def\@begineqnseries[#1]{\edef\@@eqnname{#1}}
\def\endeqnseries{\subequationnumstyle{blank}%
     \expandafter\ifnoarg\@@eqnname\then
     \else\label\@@eqnname{\puteqnformat}%
     \fi
     \aftergroup\ignorespaces}

\def\figlabel#1{%
     \if@subfigcount
          \newsubfigurenum=\next
     \else\newfigurenum=\next
          \ifblank\subfigurenumstyle\then
          \else\global\@subfigcounttrue
               \newsubfigurenum=\@ne
          \fi
     \fi
     \label{#1}{\putfigformat}\putfig{#1}%
     {\def\marginnoteformat{\tt}\marginnote{#1}}}

\let\putfig=\putlab

\def\figseriesstyle#1{\gdef\@figseriesstyle{#1}}
\def\beginfigseries{\subfigurenumstyle{\@figseriesstyle}%
     \defaultoption[]\@beginfigseries}
\def\@beginfigseries[#1]{\edef\@@figname{#1}}
\def\endfigseries{\subfigurenumstyle{blank}%
     \expandafter\ifnoarg\@@figname\then
     \else\label\@@figname{\putfigformat}%
     \fi
     \aftergroup\ignorespaces}

\def\tbllabel#1{%
     \if@subtblcount
          \newsubtablenum=\next
     \else\newtablenum=\next
          \ifblank\subtablenumstyle\then
          \else\global\@subtblcounttrue
               \newsubtablenum=\@ne
          \fi
     \fi
     \label{#1}{\puttblformat}\puttbl{#1}%
     {\def\marginnoteformat{\tt}\marginnote{#1}}}

\let\puttbl=\putlab

\def\tblseriesstyle#1{\gdef\@tblseriesstyle{#1}}
\def\begintblseries{\subtablenumstyle{\@tblseriesstyle}%
     \defaultoption[]\@begintblseries}
\def\@begintblseries[#1]{\edef\@@tblname{#1}}
\def\endtblseries{\subtablenumstyle{blank}%
     \expandafter\ifnoarg\@@tblname\then
     \else\label\@@tblname{\puttblformat}%
     \fi
     \aftergroup\ignorespaces}

\message{reference numbering,}

\newcount\referencenum \referencenum=0
\newcount\@@prerefcount \@@prerefcount=0
\newcount\@@thisref
\newcount\@@lastref
\newcount\@@loopref
\newcount\@@refseq
\newdimen\refnumindent
\let\@undefreflist=\empty

\def\referencenumstyle#1{\@setnumstyle\referencenum{#1}}

\def\referencestyle#1{\usename{@ref#1}}

\def\@refsequential{%
     \gdef\@refpredef##1{\global\advance\referencenum by\@ne
          \let\^^\=0\label{##1}{\^^\{\the\referencenum}}%
          \useafter\gdef{ref@\the\referencenum}{{##1}{\undefinedlabelformat}}}%
     \gdef\@reference##1##2{%
          \ifundefinedlabel##1\then
          \else\def\^^\####1{\global\@@thisref=####1\relax}\putlab{##1}%
               \useafter\gdef{ref@\the\@@thisref}{{##1}{##2}}%
          \fi}%
     \gdef\endputreferences{%
          \loop\ifnum\@@loopref<\referencenum
                    \advance\@@loopref by\@ne
                    \expandafter\expandafter\expandafter\@printreference
                         \csname ref@\the\@@loopref\endcsname
          \repeat
          \par}}

\def\@refpreordered{%
     \gdef\@refpredef##1{\global\advance\referencenum by\@ne
          \additemR##1\to\@undefreflist}%
     \gdef\@reference##1##2{%
          \ifundefinedlabel##1\then
          \else\global\advance\@@loopref by\@ne
               {\let\^^\=0\label{##1}{\^^\{\the\@@loopref}}}%
               \@printreference{##1}{##2}%
          \fi}
     \gdef\endputreferences{%
          \def\^^\####1{\useiflab{####1}\then
               \else\reference{####1}{\undefinedlabelformat}\fi}%
          \@undefreflist
          \par}}

\def\beginprereferences{\par
     \def\reference##1##2{\global\advance\referencenum by1\@ne
          \let\^^\=0\label{##1}{\^^\{\the\referencenum}}%
          \useafter\gdef{ref@\the\referencenum}{{##1}{##2}}}}
\def\endprereferences{\global\@@prerefcount=\the\referencenum\par}

\def\beginputreferences{\par
     \refnumindent=\z@\@@loopref=\z@
     \loop\ifnum\@@loopref<\referencenum
               \advance\@@loopref by\@ne
               \setbox\z@=\hbox{\referencenum=\@@loopref
                    \referencenumformat\enskip}%
               \ifdim\wd\z@>\refnumindent\refnumindent=\wd\z@\fi
     \repeat
     \putreferenceformat
     \@@loopref=\z@
     \loop\ifnum\@@loopref<\@@prerefcount
               \advance\@@loopref by\@ne
               \expandafter\expandafter\expandafter\@printreference
                    \csname ref@\the\@@loopref\endcsname
     \repeat
     \let\reference=\@reference}

\def\@printreference#1#2{\ifx#2\undefinedlabelformat\err@undefinedref{#1}\fi
     \noindent\ifdraft\rlap{\hskip\hsize\hskip.1in \tt#1}\fi
     \llap{\referencenum=\@@loopref\referencenumformat\enskip}#2\par}

\def\reference#1#2{{\par\refnumindent=\z@\putreferenceformat\noindent#2\par}}

\def\putref#1{\trim#1\to\@@refarg
     \expandafter\ifnoarg\@@refarg\then
          \toks@={\relax}%
     \else\@@lastref=-\@m\def\@@refsep{}\def\@more{\@nextref}%
          \toks@={\@nextref#1,,}%
     \fi\the\toks@}
\def\@nextref#1,{\trim#1\to\@@refarg
     \expandafter\ifnoarg\@@refarg\then
          \let\@more=\relax
     \else\ifundefinedlabel\@@refarg\then
               \expandafter\@refpredef\expandafter{\@@refarg}%
          \fi
          \def\^^\##1{\global\@@thisref=##1\relax}%
          \global\@@thisref=\m@ne
          \setbox\z@=\hbox{\putlab\@@refarg}%
     \fi
     \advance\@@lastref by\@ne
     \ifnum\@@lastref=\@@thisref\advance\@@refseq by\@ne\else\@@refseq=\@ne\fi
     \ifnum\@@lastref<\z@
     \else\ifnum\@@refseq<\thr@@
               \@@refsep\def\@@refsep{,}%
               \ifnum\@@lastref>\z@
                    \advance\@@lastref by\m@ne
                    {\referencenum=\@@lastref\putrefformat}%
               \else\undefinedlabelformat
               \fi
          \else\def\@@refsep{--}%
          \fi
     \fi
     \@@lastref=\@@thisref
     \@more}

\message{streaming,}

\newif\ifstreaming

\def\streamto{\defaultoption[\jobname]\@streamto}
\def\@streamto[#1]{\global\streamingtrue
     \immediate\write\sixt@@n{--> Streaming to #1.str}%
     \newwrite\streamout\immediate\openout\streamout=#1.str }

\def\streamfrom{\defaultoption[\jobname]\@streamfrom}
\def\@streamfrom[#1]{\newread\streamin\openin\streamin=#1.str
     \ifeof\streamin
          \expandafter\err@nostream\expandafter{#1.str}%
     \else\immediate\write\sixt@@n{--> Streaming from #1.str}%
          \let\@@labeldef=\gdef
          \ifstreaming
               \edef\@elc{\endlinechar=\the\endlinechar}%
               \endlinechar=\m@ne
               \loop\read\streamin to\@@scratcha
                    \ifeof\streamin
                         \streamingfalse
                    \else\toks@=\expandafter{\@@scratcha}%
                         \immediate\write\streamout{\the\toks@}%
                    \fi
                    \ifstreaming
               \repeat
               \@elc
               \input #1.str
               \streamingtrue
          \else\input #1.str
          \fi
          \let\@@labeldef=\xdef
     \fi}

\def\streamcheck{\ifstreaming
     \immediate\write\streamout{\pagenum=\the\pagenum}%
     \immediate\write\streamout{\footnotenum=\the\footnotenum}%
     \immediate\write\streamout{\referencenum=\the\referencenum}%
     \immediate\write\streamout{\chapternum=\the\chapternum}%
     \immediate\write\streamout{\sectionnum=\the\sectionnum}%
     \immediate\write\streamout{\subsectionnum=\the\subsectionnum}%
     \immediate\write\streamout{\equationnum=\the\equationnum}%
     \immediate\write\streamout{\subequationnum=\the\subequationnum}%
     \immediate\write\streamout{\figurenum=\the\figurenum}%
     \immediate\write\streamout{\subfigurenum=\the\subfigurenum}%
     \immediate\write\streamout{\tablenum=\the\tablenum}%
     \immediate\write\streamout{\subtablenum=\the\subtablenum}%
     \immediate\closeout\streamout
     \fi}


\def\err@badtypesize{%
     \errhelp={The limited availability of certain fonts requires^^J%
          that the base type size be 10pt, 12pt, or 14pt.^^J}%
     \errmessage{--> Illegal base type size}}

\def\err@badsizechange{\immediate\write\sixt@@n
     {--> Size change not allowed in math mode, ignored}}

\def\err@sizetoolarge#1{\immediate\write\sixt@@n
     {--> \noexpand#1 too big, substituting HUGE}}

\def\err@sizenotavailable#1{\immediate\write\sixt@@n
     {--> Size not available, \noexpand#1 ignored}}

\def\err@fontnotavailable#1{\immediate\write\sixt@@n
     {--> Font not available, \noexpand#1 ignored}}

\def\err@sltoit{\immediate\write\sixt@@n
     {--> Style \noexpand\sl not available, substituting \noexpand\it}%
     \it}

\def\err@bfstobf{\immediate\write\sixt@@n
     {--> Style \noexpand\bfs not available, substituting \noexpand\bf}%
     \bf}

\def\err@badgroup#1#2{%
     \errhelp={The block you have just tried to close was not the one^^J%
          most recently opened.^^J}%
     \errmessage{--> \noexpand\end{#1} doesn't match \noexpand\begin{#2}}}

\def\err@badcountervalue#1{\immediate\write\sixt@@n
     {--> Counter (#1) out of bounds}}

\def\err@extrafootnotemark{\immediate\write\sixt@@n
     {--> \noexpand\footnotemark command
          has no corresponding \noexpand\footnotetext}}

\def\err@extrafootnotetext{%
     \errhelp{You have given a \noexpand\footnotetext command without first
          specifying^^Ja \noexpand\footnotemark.^^J}%
     \errmessage{--> \noexpand\footnotetext command has no corresponding
          \noexpand\footnotemark}}

\def\err@labelredef#1{\immediate\write\sixt@@n
     {--> Label "#1" redefined}}

\def\err@badlabelmatch#1{\immediate\write\sixt@@n
     {--> Definition of label "#1" doesn't match value in \jobname.lab}}

\def\err@needlabel#1{\immediate\write\sixt@@n
     {--> Label "#1" cited before its definition}}

\def\err@undefinedlabel#1{\immediate\write\sixt@@n
     {--> Label "#1" cited but never defined}}

\def\err@undefinedeqn#1{\immediate\write\sixt@@n
     {--> Equation "#1" not defined}}

\def\err@undefinedref#1{\immediate\write\sixt@@n
     {--> Reference "#1" not defined}}

\def\err@nostream#1{%
     \errhelp={You have tried to input a stream file that doesn't exist.^^J}%
     \errmessage{--> Stream file #1 not found}}

\message{jyTeX initialization}

\everyjob{\immediate\write16{--> jyTeX version \fmtversion}%
     \edef\@@jobname{\jobname}%
     \edef\jobname{\@@jobname}%
     \settime
     \openin0=\jobname.lab
     \ifeof0
     \else\closein0
          \immediate\write16{--> Getting labels from file \jobname.lab}%
          \input\jobname.lab
     \fi}


\def\fixedskipslist{%
     \^^\{\topskip}%
     \^^\{\splittopskip}%
     \^^\{\maxdepth}%
     \^^\{\skip\topins}%
     \^^\{\skip\footins}%
     \^^\{\headskip}%
     \^^\{\footskip}}

\def\scalingskipslist{%
     \^^\{\p@renwd}%
     \^^\{\delimitershortfall}%
     \^^\{\nulldelimiterspace}%
     \^^\{\scriptspace}%
     \^^\{\jot}%
     \^^\{\normalbaselineskip}%
     \^^\{\normallineskip}%
     \^^\{\normallineskiplimit}%
     \^^\{\baselineskip}%
     \^^\{\lineskip}%
     \^^\{\lineskiplimit}%
     \^^\{\bigskipamount}%
     \^^\{\medskipamount}%
     \^^\{\smallskipamount}%
     \^^\{\parskip}%
     \^^\{\parindent}%
     \^^\{\abovedisplayskip}%
     \^^\{\belowdisplayskip}%
     \^^\{\abovedisplayshortskip}%
     \^^\{\belowdisplayshortskip}%
     \^^\{\abovechapterskip}%
     \^^\{\belowchapterskip}%
     \^^\{\abovesectionskip}%
     \^^\{\belowsectionskip}%
     \^^\{\abovesubsectionskip}%
     \^^\{\belowsubsectionskip}}


\def\twoupsetup{
     \topmargin=.75in
     \leftmargin=.5in
     \vsize=6.9in
     \hsize=4.75in
     \fullhsize=10in
     \let\draft=\relax}

\outputstyle{normal}                             

\def\marginnoteformat{\subscriptsize             
     \hsize=1in \baselinestretch=1000 \everypar={}%
     \tolerance=5000 \hbadness=5000 \parskip=0pt \parindent=0pt
     \leftskip=0pt \rightskip=0pt \raggedright}

\head={\ifdraft\normalfonts\it\hfil DRAFT\hfil   
     \llap{\number\day\ \monthword\month\ \militarytime}\else\hfil\fi}
\foot={\hfil\normalfonts\numstyle\pagenum\hfil}  

\normalbaselineskip=12pt                         
\normallineskip=0pt                              
\normallineskiplimit=0pt                         
\normalbaselines                                 

\topskip=.85\baselineskip \splittopskip=\topskip \headskip=2\baselineskip
\footskip=\headskip

\pagenumstyle{arabic}                            

\parskip=0pt                                     
\parindent=20pt                                  

\baselinestretch=1000                            


\chapterstyle{left}                              
\chapternumstyle{blank}                          
\def\chapterbreak{\newpage}                      
\abovechapterskip=0pt                            
\belowchapterskip=1.5\baselineskip               
     plus.38\baselineskip minus.38\baselineskip
\def\chapternumformat{\numstyle\chapternum.}     

\sectionstyle{left}                              
\sectionnumstyle{blank}                          
\def\sectionbreak{\vskip0pt plus4\baselineskip\penalty-100
     \vskip0pt plus-4\baselineskip}              
\abovesectionskip=1.5\baselineskip               
     plus.38\baselineskip minus.38\baselineskip
\belowsectionskip=\the\baselineskip              
     plus.25\baselineskip minus.25\baselineskip
\def\sectionnumformat{
     \ifblank\chapternumstyle\then\else\numstyle\chapternum.\fi
     \numstyle\sectionnum.}

\subsectionstyle{left}                           
\subsectionnumstyle{blank}                       
\def\subsectionbreak{\vskip0pt plus4\baselineskip\penalty-100
     \vskip0pt plus-4\baselineskip}              
\abovesubsectionskip=\the\baselineskip           
     plus.25\baselineskip minus.25\baselineskip
\belowsubsectionskip=.75\baselineskip            
     plus.19\baselineskip minus.19\baselineskip
\def\subsectionnumformat{
     \ifblank\chapternumstyle\then\else\numstyle\chapternum.\fi
     \ifblank\sectionnumstyle\then\else\numstyle\sectionnum.\fi
     \numstyle\subsectionnum.}


\footnotenumstyle{symbols}                       
\footnoteskip=0pt                                
\def\footnotenumformat{\numstyle\footnotenum}    
\def\footnoteformat{\footnotesize                
     \everypar={}\parskip=0pt \parfillskip=0pt plus1fil
     \leftskip=1em \rightskip=0pt
     \spaceskip=0pt \xspaceskip=0pt
     \def\\{\ifhmode\ifnum\lastpenalty=-10000
          \else\hfil\penalty-10000 \fi\fi\ignorespaces}}


\def\undefinedlabelformat{$\bullet$}             


\equationnumstyle{arabic}                        
\subequationnumstyle{blank}                      
\figurenumstyle{arabic}                          
\subfigurenumstyle{blank}                        
\tablenumstyle{arabic}                           
\subtablenumstyle{blank}                         

\eqnseriesstyle{alphabetic}                      
\figseriesstyle{alphabetic}                      
\tblseriesstyle{alphabetic}                      

\def\puteqnformat{\hbox{
     \ifblank\chapternumstyle\then\else\numstyle\chapternum.\fi
     \ifblank\sectionnumstyle\then\else\numstyle\sectionnum.\fi
     \ifblank\subsectionnumstyle\then\else\numstyle\subsectionnum.\fi
     \numstyle\equationnum
     \numstyle\subequationnum}}
\def\putfigformat{\hbox{
     \ifblank\chapternumstyle\then\else\numstyle\chapternum.\fi
     \ifblank\sectionnumstyle\then\else\numstyle\sectionnum.\fi
     \ifblank\subsectionnumstyle\then\else\numstyle\subsectionnum.\fi
     \numstyle\figurenum
     \numstyle\subfigurenum}}
\def\puttblformat{\hbox{
     \ifblank\chapternumstyle\then\else\numstyle\chapternum.\fi
     \ifblank\sectionnumstyle\then\else\numstyle\sectionnum.\fi
     \ifblank\subsectionnumstyle\then\else\numstyle\subsectionnum.\fi
     \numstyle\tablenum
     \numstyle\subtablenum}}


\referencestyle{sequential}                      
\referencenumstyle{arabic}                       
\def\putrefformat{\numstyle\referencenum}        
\def\referencenumformat{\numstyle\referencenum.} 
\def\putreferenceformat{
     \everypar={\hangindent=1em \hangafter=1 }%
     \def\\{\hfil\break\null\hskip-1em \ignorespaces}%
     \leftskip=\refnumindent\parindent=0pt \interlinepenalty=1000 }


\normalsize


\def\fmtversion{2.6M (June 1992)}

\catcode`\@=12

\typesize=10pt \magnification=1200 \baselineskip17truept
\footnotenumstyle{arabic} \hsize=6truein\vsize=8.5truein
\input epsf
\sectionnumstyle{blank}
\chapternumstyle{blank}
\chapternum=1
\sectionnum=1
\pagenum=0

\def\begintitle{\pagenumstyle{blank}\parindent=0pt
\begin{narrow}[0.4in]}
\def\endtitle{\end{narrow}\newpage\pagenumstyle{arabic}}


\def\beginexercise{\vskip 20truept\parindent=0pt\begin{narrow}[10
truept]}
\def\endexercise{\vskip 10truept\end{narrow}}


\def\eql#1{\eqno\eqnlabel{#1}}
\def\ref{\reference}
\def\peq{\puteqn}
\def\pref{\putref}

\def\mgn{\marginnote}
\def\bex{\begin{exercise}}
\def\eex{\end{exercise}}


\font\open=msbm10 


\def\StretchRtArr#1{{\count255=0\loop\relbar\joinrel\advance\count255 by1
\ifnum\count255<#1\repeat\rightarrow}}
\def\StretchLtArr#1{\,{\leftarrow\!\!\count255=0\loop\relbar
\joinrel\advance\count255 by1\ifnum\count255<#1\repeat}}

\def\StretchLRtArr#1{\,{\leftarrow\!\!\count255=0\loop\relbar\joinrel\advance
\count255 by1\ifnum\count255<#1\repeat\rightarrow\,\,}}

\def\mbox#1{{\leavevmode\hbox{#1}}}

\def\hspace#1{{\phantom{\mbox#1}}}
\def\oZ{\mbox{\open\char90}}


\def\Ga{\Gamma}

\def\La{\Lambda}
\def\om{\omega}

\def\ze{\zeta}

\def\De{\Delta}

\def\caO{{\cal O}}

\def\sc{{\rm sc }}

\def\zf{$\zeta$--function}
\def\zfs{$\zeta$--functions}


\def\frac#1/#2{\leavevmode\kern.1em
\raise.5ex\hbox{\the\scriptfont0 #1}\kern-.1em/\kern-.15em
\lower.25ex\hbox{\the\scriptfont0 #2}}
\def\sfrac#1/#2{\leavevmode\kern.1em
\raise.5ex\hbox{\the\scriptscriptfont0 #1}\kern-.1em/\kern-.15em
\lower.25ex\hbox{\the\scriptscriptfont0 #2}}

\def\gtorder{\mathrel{\raise.3ex\hbox{$>$}\mkern-14mu
             \lower0.6ex\hbox{$\sim$}}}
\def\ltorder{\mathrel{\raise.3ex\hbox{$<$}\mkern-14mu
             \lower0.6ex\hbox{$\sim$}}}

\def\semidirprod{\rlap{\ss C}\raise1pt\hbox{$\mkern.75mu\times$}}
\def\for{\lower6pt\hbox{$\Big|$}}
\def\fish{\kern-.25em{\phantom{abcde}\over \phantom{abcde}}\kern-.25em}


\def\boxit#1{\vbox{\hrule\hbox{\vrule\kern3pt
        \vbox{\kern3pt#1\kern3pt}\kern3pt\vrule}\hrule}}
\def\dalemb#1#2{{\vbox{\hrule height .#2pt
        \hbox{\vrule width.#2pt height#1pt \kern#1pt \vrule
                width.#2pt} \hrule height.#2pt}}}

\def\frac#1#2{{{#1}\over{#2}}}

\def\noin{\noindent}

\def\comb#1#2{{\left(#1\atop#2\right)}}

\def\eg{{\it e.g.}}

\def\cf{{\it cf }}
\def\pa{\partial}



\def\wt{\widetilde}

\def\3j#1#2#3#4#5#6{\left\lgroup\matrix{#1&#2&#3\cr#4&#5&#6\cr}
\right\rgroup}

\def\m?{\mgn{?}}

\def\pa{\partial}

\def\beq{\begin{eqnarray}}
\def\eeq{\end{eqnarray}}


\def\aop#1#2#3{{\it Ann. Phys.} {\bf {#1}} ({#2}) #3}
\def\cjp#1#2#3{{\it Can. J. Phys.} {\bf {#1}} ({#2}) #3}
\def\cmp#1#2#3{{\it Comm. Math. Phys.} {\bf {#1}} ({#2}) #3}
\def\cqg#1#2#3{{\it Class. Quant. Grav.} {\bf {#1}} ({#2}) #3}

\def\jmp#1#2#3{{\it J. Math. Phys.} {\bf {#1}} ({#2}) #3}
\def\jpa#1#2#3{{\it J. Phys.} {\bf A{#1}} ({#2}) #3}

\def\np#1#2#3{{\it Nucl. Phys.} {\bf B{#1}} ({#2}) #3}

\def\pl#1#2#3{{\it Phys. Lett.} {\bf {#1}} ({#2}) #3}

\def\prB#1#2#3{{\it Phys. Rev.} {\bf B{#1}} ({#2}) #3}
\def\prD#1#2#3{{\it Phys. Rev.} {\bf D{#1}} ({#2}) #3}
\def\prl#1#2#3{{\it Phys. Rev. Lett.} {\bf #1} ({#2}) #3}

\def\am#1#2#3{{\it Acta Mathematica} {\bf {#1}} ({#2}) #3}
\def\aim#1#2#3{{\it Adv. in Math.} {\bf {#1}} ({#2}) #3}
\def\ajm#1#2#3{{\it Am. J. Math.} {\bf {#1}} ({#2}) #3}

\def\jpamt#1#2#3{{\it J. Phys.A:Math.Theor.} {\bf{#1}} ({#2}) #3}
\def\jram#1#2#3{{\it J. f. reine u. Angew. Math.} {\bf {#1}} ({#2}) #3}

\def\ma#1#2#3{{\it Math. Ann.} {\bf {#1}} ({#2}) #3}

\def\mz#1#2#3{{\it Math. Zeit.} {\bf {#1}} ({#2}) #3}

\def\plb#1#2#3{{\it Phys. Letts.} {\bf {B#1}} ({#2}) #3}

\def\qjm#1#2#3{{\it Quart. J. Math.} {\bf {#1}} ({#2}) #3}

\def\rmjm#1#2#3{{\it Rocky Mountain J. Math.} {\bf {#1}} ({#2}) #3}

\def\tams#1#2#3{{\it Trans.Am.Math.Soc.} {\bf {#1}} ({#2}) #3}

\begin{title}
\vglue 0.5truein
\vskip15truept
\centertext {\Bigfonts \bf Computation of the derivative } \vskip7truept
\vskip10truept\centertext{\Bigfonts \bf of the Hurwitz \zf\ and}
 \vskip7truept
\vskip10truept\centertext{\Bigfonts \bf the higher Kinkelin constants}
 \vskip 20truept
\centertext{J.S.Dowker\footnote{dowker@man.ac.uk;  dowkeruk@yahoo.co.uk}} \vskip
7truept \centertext{\it Theory Group,} \centertext{\it School of Physics and Astronomy,}
\centertext{\it The University of Manchester,} \centertext{\it Manchester, England} \vskip
7truept \centertext{}

\vskip 7truept \vskip40truept
\begin{narrow}
I use the numerical values of the generalised Glaisher-Kinkerlin-\break Bendersky (GKB)
constants to give numerical values for the derivatives of the Hurwitz \zf\ at negative
integers, rather than the other way round. I point out that both Glaisher's numerical
approach and Bendersky's recursion for the generalised gamma function were anticipated
by Jeffery in 1862 who gave the value of the second constant as an example. I therefore
propose that GKB become GKBJ.

\end{narrow}
\vskip 5truept
\vskip 60truept
\vfil
\end{title}
\pagenum=0
\newpage

\section{\bf 1. Introduction}

In many calculations of the effective action, long combinations of derivatives of  Hurwitz
\zfs\ occur and often it can be more useful to give their numerical values. Although this is
a fairly standard topic, some further discussion might be of interest. One could of course
use a built--in algorithm in a suitable CAS, but if this is not available, a direct analytical
calculation is necessary, preferably not involving quadrature. Suitable forms have been
found by Elizalde, [\pref{Elizalde}], some time ago. My calculation will incidentally
produce the same (asymptotic) expression but by a different route.

\section{\bf2. The basics}
I approach the formula in a possibly circular  way allowing me to introduce some other
notions and, later, a few historical remarks. I start with the basic relation between
$\ze_H'(-k,w)$ and $\ze_H'(0,1)$, the Riemann \zf, with (initially) $w$ an integer,
  $$
  \ze'(-k,w)-\ze'(-k)=\log \Ga_k(w)\,,
  \eql{bas}
  $$
dropping the `H'. Here $\Ga_k(w)$ is Bendersky's generalised gamma function defined
(initially) by the sum
  $$
  \log\Ga_k(w+1)=\sum_{m=1}^w m^k\,\log m\,,\quad w\in\oZ\,.
  \eql{gengam}
  $$

There are two attitudes one can take towards (\peq{bas}). These are, in essence, whether
it is to be read from right to left or left to right. In the first case, which is the one I mostly
adopt here, the right--hand side, defined as the explicit sum (\peq{gengam}), is analysed
independently using finite calculus summation techniques leading to calculable quantities
in which $w$ can be considered continuous. This is the procedure of Bendersky,
[\pref{bend}], and Jeffery, [\pref{Jeffery}]. (They do not consider the relation with the
\zf.) Having an expression for the generalised gamma function then allows information
about the \zf\ to be extracted from (\peq{bas}).

On the other hand  if the \zf\ is assumed known  then $\Ga_k$ can be taken as {\it
defined} by (\peq{bas}) and its properties determined therefrom. I will not consider this
point of view much here but it has the advantages of rapidity and elegance. An example
will appear later.

\begin{ignore} If the left--hand side of (\peq{bas}) is known for all $w$, it can be taken
as the definition of $\Ga_k(w)$. This is the point of view of Kurokawa and Ochiai,
[\pref{KandO}], who are unaware of Bendersky's work, but I would prefer to run the
equation the other way and {\it calculate} the left--hand side.

Equation (\peq{bas}) suggests a generalisation which is obtained by replacing Hurwitz by
Barnes \zfs. This in pursued by Kurokawa and Ochiai.
\end{ignore}

I now invoke, without worrying too much about where it comes from, Adamchik's formula
for $\ze'(-k)$, [\pref{Adamchik1}],
  $$
  \ze'(-k)={H_k\,B_{k+1}\over{k+1}}-\log A_k\,,
  \eql{adam}
  $$
where $A_k$ are Bendersky's generalisation of the Glaisher--Kinkelin constants defined by
the asymptotic behaviour,
  $$
 \log A_k\equiv L_k= \lim_{w\to\infty} \log\Ga_k(w+1)\big|_{w{\rm -independant} }\,.
 \eql{limit}
  $$
$H_k$ is a harmonic number.

One writes,
  $$
  \log\Ga_k(w+1)=L_k+\La_k(w+1)
  \eql{bend1}
  $$
where $\La_k(w)$ is a function with the appropriate limiting behaviour. Bendersky derives
an explicit asymptotic (divergent)  series from the Euler--Maclaurin summation formula
applied to $\log\Ga_k(w+1)$. It is not necessary to write it out yet.

Combining (\peq{bas}), (\peq{adam}) and (\peq{bend1}) trivially yields the required
computable formula,
  $$\eqalign{
  \ze'(-k,w)&={H_k\,B_{k+1}\over{k+1}}-L_k+\log\Ga_k(w)\cr
  &={H_k\,B_{k+1}\over{k+1}}+\La_k(w)\,,
  }
  \eql{ashzet}
  $$
which gives the asymptotic expansion obtained by Elizalde directly from integral forms of
the \zf.

I look upon (\peq{adam}) as a means of finding $\ze'(-k)$ from the $A_k$, rather than the
other way round, as is more usual. As a check of the numbers, if $k$ is even $\ze'(-k)$ can
be transformed into $\ze(2k+1)$ the values of which are readily available, to high
accuracy.\footnote{Or say Wolfram Alpha can be employed.}

Numerically one can treat the two lines in (\peq{ashzet}) separately. In the top line $L_k$
is to be calculated by trial from (\peq{limit}) using the known asymptotic form of $\La$.
(The method used by Jeffery, Glaisher and Bendersky). Less accurately, one can just
substitute the asymptotic form directly into the second line. Of course, using an
asymptotic form is not the best numerical procedure, but an accuracy of more than 10
places is easily attained this way. (The isolated case of $w=1$ is not possible in the
second approach.)

The more accurate computation follows from the first line in (\peq{ashzet}) where $L_k$
is found from (\peq{limit}) using the asymptotic form for $\La_k(w+1)$ with a suitable
choice for $w$. $\log\Ga_k(w)$ is determined `exactly' by (\peq{gengam}). This method is
therefore restricted to $w$ integral.

If $w$ is not integral, one has to use the second line in (\peq{ashzet}). If $w$ is small a
direct application is not possible but one can translate $w$ to a large enough value by
adding an integer and then employing the fundamental property of the generalised Gamma
functions, [\pref{bend}], valid for any $x$,
  $$
\Ga_k(x+1)=x^{x^k}\,\Ga_k(x)\,,
\eql{fundprop}
  $$
several times, if necessary

\section{\bf3. The details}

I will now fill out the previous discussion by outlining Bendersky's approach. I do this
because his paper is not widely recognised. To begin with, it is helpful to write down the
explicit expression for the asymptotic series, $\La_k$, as given by Bendersky. One
particular form of this is, and I write it out exactly as in [\pref{bend}] p.276, (with a few
misprints corrected),
  $$\eqalign{
  \La_{k+1}&(x+1)={x^{k+2}\over k+2}\,\log x-{x^{k+2}\over(k+2)^2}
  +{1\over2}\,x^{k+1}\,\log x+\cr
  &+(k+1)!\sum_{r=1}^{k-1}{B_{r+1}\over(r+1)!}{x^{k+1-r}\over(k+1-r)!}
  \bigg(\log x +{1\over k+1}+{1\over k}+\ldots+{1\over k+2-r}\bigg)+\cr
  &+x\,B_{k+1}\bigg(\log x+{1\over k+1}+{1\over k} +\ldots+{1\over3}+{1\over2}
  \bigg)+{B_{k+2}\over k+2}\,\log x+\cr
  &+(k+1)!\sum_{s=1}^\infty (-1)^s{B_{k+1+s}\over (k+1+s)!}{(s-2)!\over x^{s-1}}\,,
  }
  \eql{lambda}
  $$
which is easily coded.

Bendersky derives this, at general $k$, from an Euler--Maclaurin summation, after some
amalgamation of terms. The details are of no immediate concern.

It is always illuminating to have some specific examples before one's eyes and I give the
following,
  $$\eqalign{
  \La_2(x+1)&={x(x+1)(2x+1)\over 1.2.3}\,\log x-{x^3\over9}+{x\over12}-{1\over360\,x}
  +{1\over7560 x^3}-\ldots\cr
  \La_1(x+1)&=\bigg({x(x+1)\over 1.2}+{1\over12}\bigg)\log x-{x^2\over4}+{1\over 720x^2}-
  {1\over5040x^4}+\ldots\cr
  \La_0(x+1)&=\bigg(x+{1\over2}\bigg)\log x-x+{1\over12x}-{1\over360x^3}+
  {1\over1260x^5}-\ldots\cr
  \La_{-1}(x+1)&=\log x+{1\over2x}-{1\over12x^2}+{1\over120 x^4}-\ldots\cr
  \La_{-2}(x+1)&=-{1\over x}+{1\over2x^2}-{1\over6x^3}+{1\over30x^5}-\ldots\,.
  }
  \eql{lambdaex}
  $$

I have included two `lower' expansions which can be taken as part of the set. Explicitly,
   $$\eqalign{
   \La_{-1}(x+1)&={1\over1}+{1\over2}+\ldots+{1\over x}-S_1\cr
   \La_{-2}(x+1)&={1\over1^2}+{1\over2^2}+\ldots+{1\over x^2}-S_2\cr
    }
   $$
where $S_1$ and $S_2$ are the corresponding `Glaisher--Kinkelin--Bendersky (GKB)'
constants and have the (known) values,
   $$
    S_1=-\Ga'(1)\,,\quad S_2=\zeta(2)\,.
   $$

We recognise in (\peq{lambdaex}) some standard asymptotic (Stirling) series, \eg\ for,
  $$
\Lambda_0(x+1)=\log{\Ga(x+1)\over \sqrt{2\pi}}\,,\quad \Lambda_{-1}(x+1)
={d\over dx}\log\Ga(x+1)=\psi(x+1)\,.
  $$
and so on downwards to give the polygamma functions, up to a factor.

As is visually clear, one can run up and down the right--hand sides in (\peq{lambdaex}) by
integration and differentiation. The exact relationship is given by Bendersky who derived it
by brute force from (\peq{lambda}). He found,
  $$
\La_{k+1}(x+1)=(k+1)\int \La_k(x+1)\,dx+{1\over k+1}\phi_{k+1}(x+1)+x\,H_k\,B_{k+1}\,,
\eql{recurs}
  $$
where $\phi_n(x)$ is a Bernoulli polynomial given in terms of the more usual polynomials,
\eg\ N\"orlund, [\pref{Norlund}], $B_n(x)$, by,

  $$
    \phi_n(x)={1\over n+1}\big(B_{n+1}(x)-B_{n+1}\big)\,,
  $$
which equals the sum of the $n$th powers of the first $x-1$ integers.

Constants of integration can be considered to be absorbed into the GKB constants, $L_k$,
which are found, in each case, by trial of varying $x$.

It is important to note that, although derived using an asymptotic series, the recursion,
(\peq{recurs}), is valid generally.

Another derivation of (\peq{recurs}) is given later.

The actual numerical values of the GKB constants, $L_k$, can be determined from the
definition, (\peq{bend1})
  $$
    L_k=\log\Ga(w+1)-\La_k(w+1)
  $$
using (\peq{gengam}) and a suitably chosen value for the {\it integer}, $w$. This choice
is linked to the necessary truncation of the infinite series in (\peq{lambda}) and leads to
an accuracy of at least 29 places, \eg\ for $w=100$ and 20 terms of the sum retained.

Bendersky also developes recursion relations for the $L_k$ which involve just convergent
series.

In order to find a more precise definition of $\log\Ga_k$ than the asymptotic series, one
starts from the observation that $\Ga_0(x+1)=\Ga(x+1)$, in terms of the ordinary Gamma
function and takes this as the starting point of an upwards recursion based on
(\peq{recurs}). To take the simplest case, $\log\Ga_1(x+1)$ has to involve the integral of
$\log(\Ga(x+1)/\sqrt{2\pi})$ and now, for exactness, a definite integral is used,
  $$
  \log\Ga_1(x+1)=\int_0^xdx\,\log{\Ga(x+1)\over\sqrt{2\pi}}+X(x)\,,
  $$
where $X$ is to be found. Because (\peq{recurs}) has to hold up to a constant, $C$,
  $$
  X=\phi_1(x+1)-C={1\over2}x(x+1)-C\,,
  $$
which is determined by setting $x=0$, so that,\footnote{ The particular values
$\log\Ga_k(1)=\log\Ga_k(2)=0$ have been used as boundary conditions.}
  $$\eqalign{
  C&=X(0)\cr
  &=0\,,
  }
  $$
and, finally,
  $$\eqalign{
  \log\Ga_1(x+1)&=\int_0^xdx\,\log\Ga(x+1)+\phi_1(x+1)-x\,L_0\cr
  &=\int_0^xdx\,\log\Ga(x+1)+{x(x+1)\over2}-x\log\sqrt{2\pi}
  }
  \eql{Raabe}
  $$
defining $\Ga_1$ in terms of the standard $\Ga$. (Setting $x=1$ yields a known result.)

Barnes, [\pref{BarnesGf}] p.281, refers to this equation (or its equivalent) as {\it
Alexeiewsky's theorem}.

 Just to
see how things fit together, I look at the next iteration according to (\peq{recurs}),
  $$\eqalign{
  \log\Ga_2(x+1)=2\int_0^x&dx\,\log\Ga_1(x+1)+{1\over2}\phi_2(x+1)+{1\over6}x
  -2xL_1\cr
=2\int_0^x&dx\int_0^x dx'\log\Ga(x'+1)+2\int_0^x dx\big(\phi_1(x+1)-xL_0\big)+\cr
&+{1\over2}\phi_2(x+1)+{1\over6}x
  -2xL_1\cr
=2\int_0^x&dx\int_0^x dx'\log\Ga(x'+1)+{3\over2}\phi_2(x+1)-2xL_1-x^2L_0\,.
  }
  $$

A continuation of this process produces Bendersky's important general solution,
  $$
  \log \Ga_k(x+1)=k!\,I_k(x)+H_k\,\phi_k(x+1)-\psi_k(x)\,,
  \eql{defgengam}
  $$
with the definitions
  $$\eqalign{
  I_k(x)&=\int_0^x dx'\, I_{k-1}(x')\,,\quad I_0(x)=\log\Ga(x+1)\cr
  \psi_k(x)&=\sum_{r=0}^{k-1}\comb kr \,L_r\,x^{k-r}\,.
  }
  \eql{fdefs}
  $$

Bendersky takes the convergent (\peq{defgengam}) as the definition of the generalised
Gamma functions for {\it all} $x$ and any integer $k$ and uses it systematically to obtain
their basic properties such as (\peq{fundprop}), encountered earlier. It would be out of
place to give a summary here.

\section{\bf 4. Earlier history}

The quantities $\log\Ga_k(x)$  in (\peq{gengam}) were first considered by Kinkelin who
concentrated on just $\Ga_1(x)$ in [\pref{kink}],\footnote{ Kinkelin's work is dated 1856,
but was published in 1860. I have not been able to access his earlier paper.} and derived
its properties {\it ab initio} from the expression in terms of the ordinary Gamma function,
  $$
  \log\Ga_1(x)=\int_0^x\log \Ga(t)\,dt+{x(x-1)\over2}-{1\over2}x\,\log 2\pi
  \eql{kinkdef}
  $$
which he obtained using Raabe's formula and the fact that it reduces to the summation,
(\peq{gengam}), for $x$ integral.

He computes what is essentially the constant $L_1$, not by the trial method but from
derived convergent series. These series can also be found in Bendersky.

Glaisher, [\pref{Glaisher}], gave a numerical treatment of  $L_1$ using the trial method.
In this, however, he was forestalled by the earlier work of Jeffery \footnote{ A very brief
biography can be found in [\pref{Dow}].} in a known paper [\pref{Jeffery2}], written in
1860. In a second, neglected paper, [\pref{Jeffery}], written in 1862, Jeffery gives what
seems to be the earliest calculation of the second GKB constant, $L_2$. Indeed he
describes, {\it inter alia}, the construction of the complete general recursion system,
(\peq{recurs}), and I have taken equation (\peq{lambdaex}) from his paper. It is clear he
could have computed as many of the $L_k$ as desired and so I propose that GKB be
extended to GKBJ.

It is worth noting, with Jeffery, that the coefficient of $\log x$ in $\La_k(x+1)$ equals
$\sum_{i=1}^x i^k$ or $B_{1+k}(x)/(1+k)$, although he does not seem to prove this in
generality. Jeffery's approach is worth re--exposure. It does not overtly use the
Euler--Maclaurin formula.

\section{\bf5. Jeffery's treatment}

For shortness, I define
  $$
  v_k(x)=\log\Ga_k(x+1)
  \eql{def}
  $$

From the original expression (\peq{gengam}) algebra gives,
  $$
  \De_x\,v_k(x)=(x+1)^k\,\log(x+1)\,,
  $$
which defines a system of equations.

Extending $x$ into the reals gives
  $$\eqalign{
  D\De\,v_k(x)&=\De\,D\,v_k(x)=k(x+1)^{k-1}\log(x+1)+(x+1)^{k-1}\cr
  &=k\De\,v_{k-1}(x)+(x+1)^{k-1}\,,
  }
  $$
so that, by summation,
  $$\eqalign{
  Dv_k(x)&=k\De^{-1}\De\,v_{k-1}(x)+\De^{-1}(x+1)^{k-1}+\varpi_1(k)\cr
  &=kv_{k-1}(x)+\phi_{k-1}(x+1)+\varpi(k)\cr
  }
  \eql{summ}
  $$

The final periodic constant $\varpi(k)$ is given by setting $x=0$ whuch yields
  $$\eqalign{
  \varpi(k)&=D\,v_k(0)\,,\quad k=2,3,\ldots\cr
  \varpi(1)&=D\,v_1(0)-1
  }
  \eql{pconsts}
  $$
since $v_k(0)=0=\phi_k(1)$. I will use (\peq{pconsts})  to find $v_k(0)$ rather than the
$\varpi(k)$.

Jeffery, [\pref{Jeffery2}], gives only the numerical value $\varpi(2)=-0.2475089541\ldots$
which is $1/4-2\,L_1$, in terms of the GKBJ constant for $k=1$. This can be shown in the
following somewhat lengthy way. How Jeffery works out the value is not clear to me.

\begin{ignore}
    *********

CHECKING:

Can we get \S42 from \S44?

Check Jeffery \S44.

(1) What is $v_x$?. Work backwards.
  $$
 2 v_x=Du_x-{x(x+1)\over2}+{.2475}
  $$

Numerical fact
  $$
  .247509=-{1\over4}+ 2\,L_1
  $$

To find the functional form, the asymptotic expressions are adequate. (IS THIS TRUE?)

First, go back to Bendersky's fundamental identity equn.(14) and check by machine. Seems
OK. Differentiating (14) to compare with Jeffery,

He says
  $$\eqalign{
2\La_1(x+1)&=\La_2'(x+1)-{1\over2}\phi_2'(x+1)-{1\over4}\cr
&=\La_2'(x+1)-\phi_1(x+1)-{1\over12}-{1\over6}\cr
&=\La_2'(x+1)-{x(x+1)\over2}-{1\over12}-{1\over6}\cr
&=\La_2'(x+1)-{x(x+1)\over2}-{1\over4}\,.\cr
}
  $$

Comparing  with Jeffery \S44, above
  $$\eqalign{
 2 v_x&=Du_x-{x(x+1)\over2}-{1\over4}+2L_1\cr
 }
  $$
we have that
  $$
  2v_x-2L_1=2\La_1(x+1)
  $$
or
  $$
v_x=L_1+\La_1(x+1)=\log\Ga_1(x+1)\equiv v_1(x)
  $$

Also the derivative of $\La_2(x+1)$ is equal to that of $\log\Ga_2(x+1)$ and so one can
write
  $$
  2v_1(x)=Dv_2(x)-{x(x+1)\over2}-{1\over4}+2L_1
  $$
which is a way of expressing Jeffery's relation. One still has to show that the constant,
$\varpi(k)$, is $1/4-2L_1$.

END OF CHECKING

 ******************************

QUERY:

Can we just derive (\peq{kinkdef}) using Jeffery?

Logically we can proceed as follows.

Define, as in Bendersky,
  $$
  v_k(x)=L_k+\La_k(x+1)
  $$
where $L_k$ is the (ultimately determined) constant in the limiting form of $v_k$ and
$\La_k(x+1)$ is the remainder, not necessarily asymptotic.

\end{ignore}

Integrating (\peq{summ}),
  $$
  v_k(x)-C(k)=k\int^x _0dx\,v_{k-1}(x)+{1\over k}\big(\phi_k(x+1)-B_{k} \,x\big)
  +\varpi(k)\,x
  \eql{reln}
  $$
where $C(k)=0$, again on setting $x=0$.

As a trivial check, set $k=1$. Then
$$\eqalign{
  v_1(x)&=\int^x _0dx\,v_{0}(x)+\big(\phi_1(x+1)-B_{1} \,x\big)+\varpi(1)\,x\cr
  &=\int^x _0dx\,v_{0}(x)+{x(x+1)\over2}+(\varpi(1)+{1\over2})\,x\,.
  }
  $$

Setting $x$ to zero produces nothing while $x=1$ gives, using Raabe's formula,
  $$\eqalign{
  0&=\int^1 _0dx\,v_0(x)+{3\over2}+\varpi(1)\cr
  &=\log\sqrt{2\pi}+\varpi(1)+{1\over2}
  }
  $$therefore the relation to the GKBJ constant, $L_0$, is
  $$
\varpi(1)=-L_0-{1\over2}=-\log\sqrt{2\pi}-{1\over2}\,.
  $$

Then
  $$\eqalign{
  v_1(x)&=\int^x _0dx\,v_{0}(x)+\big(\phi_1(x+1)-B_{1} \,x\big)+\varpi(1)\,x\cr
  &=\int^x _0dx\,v_{0}(x)+{x(x+1)\over2}-\log\sqrt{2\pi}\,x\,,
  }
  \eql{vee1}
  $$
which is the same as (\peq{Raabe}), or (\peq{kinkdef}).

Incidentally, from (\peq{pconsts}) , one sees for that the series for $v_1$ in ascending
powers of $x$ begins with the term $\big(1/2-\log\sqrt{2\pi}\big)\,x$.

Relation (\peq{reln}) can be iterated. Thus
  $$\eqalign{
  v_k(x)&=k\int^x _0dx'\,\bigg[(k-1)\int^{x'} _0dx''\,v_{k-2}(x'')+{1\over k-1}
  \big(\phi_{k-1}(x'+1)-B_{k-1} \,x'\big)+\cr
  &\hspace{*********}+\varpi(k-1)\,x'\bigg]+
  {1\over k}\big(\phi_k(x+1)-B_{k} \,x\big)+\varpi(k)\,x\cr
  &=k(k-1)\int^x _0dx'\,\int^{x'} _0dx''\,v_{k-2}(x'')
  +{1\over k-1}\big(\phi_k(x+1)-B_k\,x\big)\cr
 &\hspace{***}-{k\over2(k-1)}B_{k-1}\,x^2+{k\over2}\varpi(k-1)\,x^2
 +{1\over k}\big(\phi_k(x+1)-B_{k} \,x\big)+\varpi(k)\,x\,.
  }
  \eql{iter}
  $$

The simplest, non--trivial,  case is $k=2$, when,
  $$\eqalign{
 v_2(x) &=2\int^x _0dx'\,\int^{x'} _0dx''\,v_{0}(x'')+\big(\phi_2(x+1)-B_2\,x\big)\cr
 &\hspace{***}-B_{1}\,x^2+\varpi(1)\,x^2
 +{1\over 2}\big(\phi_2(x+1)-B_{2} \,x\big)+\varpi(2)\,x\,.
  }
  $$
Setting $x=1$ enables an expression for the summation constant, $\varpi(2)$, to be
found,
$$\eqalign{
 0&=2\int^1 _0dx'\,\int^{x'} _0dx''\,v_{0}(x'')+\big(\phi_2(2)-B_2\big)\cr
 &\hspace{***}-B_{1}+\varpi(1)
 +{1\over 2}\big(\phi_2(2)-B_{2} \big)+\varpi(2)\cr
&=2\int^1 _0dx\,(1-x)\,\log\Ga(x+1)+\big(\phi_2(2)-B_2\big)\cr
 &\hspace{***}-B_{1}+\varpi(1)
 +{1\over 2}\big(\phi_2(2)-B_{2} \big)+\varpi(2)\cr
&=2\int^1 _0dx\,(1-x)\,\log\Ga(x+1)+{7\over4}
+\varpi(1)+\varpi(2)\,.\cr
 }
 \eql{vee2})
 $$

 Rewrite this in terms of $\Ga(x)$,
 $$\eqalign{
&=-2\int^1 _0dx\,x\,\log\Ga(x)-{3\over2}+2\log\sqrt{2\pi}
+{7\over4}+\varpi(1)+\varpi(2)\cr
 &=-2\int^1 _0dx\,x\,\log\Ga(x)-{3\over4}-\varpi(1)
 +\varpi(2)\,.\cr
  }
  \eql{vee3}
  $$
In order now to relate the constant of integration, $\varpi(2)$, to the GKBJ asymptotic
constant, $L_1$, I use Kinkelin's general approach.

The derivation is based on the classic expression, or definition,  of the gamma function as
an infinite product. This gives,
  $$
\log\Ga(x)=\log w!+(x-1)\log w -\sum_{i=0}^{w-1}\log(x+i)\,,\quad w\to\infty\,.
\eql{iprod2}
  $$
Also needed is the asymptotic limit, assumed known, of the factorial (Stirling),
  $$
  \log w!\sim-w+w\log w +{1\over2}\log w +L_0\,.
  \eql{stir}
  $$

The relevant integral in (\peq{vee3}) is
  $$
    2\int^1 _0dx\,x\,\log\Ga(x)=\log w!-{1\over3}\log w-2\sum_{i=0}^{w-1}
    \int_0^1dx\,x\,\log(x+i)\,,
    \eql{rel}
   $$
and so one needs,
  $$
2\sum_{i=0}^{w-1}\int_0^1dx\,x\,\log(x+i)=\sum_{i=0}^{w-1}
\bigg(i^2\big(\log i-\log(i+1)\big)+\log(i+1)+i-{1\over2}\bigg)\,.
  $$
Rorganising,
  $$\eqalign{
  \sum_{i=0}^{w-1}\bigg(i^2\big(\log i-\log(i+1)\big)
  +\log(i+1)\bigg)&=\sum_{i=1}^{w-1}i^2\log i-
  \sum_{i=1}^{w}(i^2-2i+1-1)\log i\cr
  &=\sum_{i=1}^{w}2i\log i-w^2\log w\,,
  }
  $$
giving, for this piece,
  $$
  2\sum_{i=0}^{w-1}\int_0^1dx\,x\,\log(x+i)
  =\sum_{i=1}^{w}2i\log i-w^2\log w+{w^2\over2}-w\,.
  \eql{sumi}
  $$

Combining this expression with (\peq{rel}) and (\peq{stir}) I find,
  $$\eqalign{
  2\int^1 _0dx\,x\,&\log\Ga(x)=\cr
&w\log w+{\log w\over6}-{1\over2}-\varpi(1)-
  \sum_{i=1}^{w}2i\log i+w^2\log w-{w^2\over2}\,.
  }
  $$

Then (\peq{vee3}) becomes
  $$\eqalign{
  -{3\over4}&-\varpi(1)+\varpi(2)=\cr
  &w\log w+{\log w\over6}-{1\over2}-\varpi(1)-
  \sum_{i=1}^{w}2i\log i+w^2\log w-{w^2\over2}
  }
  \eql{rel3}
  $$
or
  $$
  \sum_{i=1}^{w}i\log i=\bigg({1\over12}+{1\over2}\big(w+w^2\big)\bigg)\log w
  +{1\over8}-{w^2\over4}-{1\over2}\varpi(2)
  $$
showing that the constant, $L_1$ is, by definition,
  $$
  L_1={1\over8}-{1\over2}\varpi(2)\,.
  $$

The initial conditions (\peq{pconsts}) also then show that the power series expansion of
$v_2(x)\equiv\log\Ga_2(x+1) $ starts with $(1/4-2L_1)\,x$ as stated, numerically, by
Jeffery, [\pref{Jeffery}] \S43.

Kinkelin derives a constant, $\varpi$, as a convergent series, from multiplication properties
of $\Ga_1$ and then shows that this is the constant in the asymptotic form. It can be found
from the information already obtained if I assume Kinkelin's equn.($22^a$) which is,
  $$\eqalign{
\log\varpi&=-{1\over6}+\log\sqrt{2\pi}-2\int^1 _0dx\,x\log\Ga(x)\cr
&=-{2\over3}-\varpi(1)+{3\over4}+\varpi(1)-\varpi(2)\cr &={1\over12}-\varpi(2)
=2L_1-{1\over6}\sim0.33084228740
}
$$
and agrees with Kinkelin's value, obtained by series summation.

Such a case by case approach, directly applied, is not algebraically efficient for the  higher
iterations. As two final examples I find that,
  $$
 \varpi(3)=-3L_2\sim-0.091345371176\,,\quad \varpi(4)=-{5\over72}-4L_3\sim
 0.013180972097\,.
  $$

\section{\bf6. The general solution. Some integrals}
I briefly return to Bendersky's general solution, (\peq{defgengam}), which he takes as the
definition of the function. The ingredients are in (\peq{fdefs}), the main part of which is the
nested integral $I_k$ which can be reduced to a single (fractional) integral of $\log\Ga$,
  $$
  I_k(x)={1\over (k-1)!}\int_0^x dt\,(x-t)^{k-1}\,\log\Ga(t+1)\,.
  \eql{ik}
  $$

By giving $x$ particular values, expressions for this integral can be found. For example if
$x=1$ the left--hand side of (\peq{defgengam}) vanishes and so,
  $$\eqalign{
  k\int_0^1 &dt\,(1-t)^{k-1}\log \Ga(t+1)=\psi_k(1)-H_k\,\phi_k(2)\cr
  &=\sum_{r=0}^{k-1}\comb kr\,L_r-H_k\bigg[{1\over2}
  +{1\over k+1}\bigg(1+\sum_{r=2}^k \comb{k+1}r\,B_r\bigg)\bigg]\cr
  &=\sum_{r=0}^{k-1}\comb kr\,\bigg({H_rB_{r+1}\over r+1}-\ze'(-r)\bigg)
  -H_k\bigg[{1\over2}
  +{1\over k+1}\bigg(1+\sum_{r=2}^k \comb{k+1}r\,B_r\bigg)\bigg]\,.
  }
  \eql{gint}
  $$

As a check take  $k=2$. Then simple arithmetic yields,
  $$\eqalign{
2\int_0^1dt\,(1-t)\log \Ga(t)&=-\ze'(0)-2\ze'(-1)+{1\over6}\,,\cr
}
  $$
which agrees with the value given in [\pref{DandKi}] p.674. In fact equn.(109) of
[\pref{DandKi}] contains (\peq{gint}). Other integration ranges can be accommodated.

I note that $I_k$, (\peq{ik}), is directly equivalent to the polygamma function of negative
order, $\psi^{(-n)}$, discussed by, \eg, Adamchik, [\pref{Adamchik1}]. See also Espinosa
and Moll, [\pref{EandM}].

More precisely, just changing $\Ga(t+1)$ to $\Ga(t)$,
  $$\eqalign{
  I_k(x)=&{1\over (k-1)!}\int_0^x dt\,(x-t)^{k-1}\,\log\Ga(t)
  +{1\over (k-1)!}\int_0^x dt\,(x-t)^{k-1}\,\log t \cr
   &=\psi^{(-k-1)}(x)+{1\over k!}\int_0^x ((x-t)^{k}-x^k)\,t^{-1}dt +x^k\log x\cr
   &=\psi^{(-k-1)}(x)+{1\over k!}
   x^k\bigg(\sum_{n=1}^k{(-1)^n\over n}\comb kn +\log x\bigg)\,.\cr
  }
  $$

Furthermore, Proposition 2 in [\pref{Adamchik1}] for the polygamma function  is
equivalent to the general solution (\peq{defgengam}) together with the basic relation,
(\peq{bas}).

\section{\bf7. Comments}

Unless I am missing something, Jeffery's method of section 5 does not easily give the {\it
general} form for $v_k$. Furthermore, in the details of the algebra, cancellations occur
which make $\varpi(k)$ a function of $L_{k-1}$ only, which ought to be derivable directly.

By contrast, Bendersky's approach is a top down one in that a {\it general} recursion is
found, somewhat out of the blue, for $\La_k=v_k-L_k$. This leads, essentially by
induction, to the form of $v_k$ for all $k$,  (\peq{defgengam}), which is then taken as
the definition of $v_k(x)$ for any $x$. Equation (\peq{defgengam}) straightaway shows
that,
  $$
   Dv_k(0)=H_kB_k-kL_{k-1}
  $$
as confirmed above in particular cases. A derivation of this expression that does not
depend on the explicit form of the general solution for $v_k$ would be welcomed.
\section{\bf8. Alternative approach}
I mentioned that looking at the basic equation (\peq{bas}) as the {\it definition} of
$\Ga_k(x)$ allows a rapid derivation of its properties and relations. As an example, I now
derive the recursion (\peq{recurs}) which follows almost immediately from the standard
formula, used in [\pref{DandKi}] for a related purpose (see also [\pref{DoandKi}] for a
generalisation),
     $$
  {\pa\over\pa w}\ze(s,w)=-s\,\ze(s+1,w)\,,
  \eql{diffw}
  $$
so that, differentiating with respect to $s$,
  $$
  {\pa\over\pa w}\ze'(s,w)=-\ze(s+1,w)-s\,\ze'(s+1,w)\,.
  $$
Then, setting $s=-k$ with $k$ a positive integer, and integrating, gives the result,
  $$
  (k+1)\int_0^x \bigg(\ze'(-k,t)-{\ze(-k,t)\over k +1}\bigg)\,dt=\ze'(-k-1,x)-\ze'(-k-1)\,.
  $$
where I have used $\ze'(-k,0)=\ze'(-k)$, [\pref{DandKi}], Appx.C.

This recursion is equivalent to Bendersky's, (\peq{recurs}), as a few lines of algebra,
utilising (\peq{adam}),  confirms. Kurokawa and Oshiai, [\pref{KandO}], unaware of the
work of Bendersky, use this approach and also generalise it to Barnes \zfs. I note that
their definition of the generalised gamma function (`higher depth' gamma function) does
not include the (constant) second term on the left--hand side of (\peq{bas}).

As a minor point, one of the integrations can be done using the well known relation,
  $$
 \ze(-k-1,x)-\ze(-k-1)=(k+1)\int_0^x dt\,\ze(-k,t)\,,\quad k>0\,,
  $$
which follows trivially from (\peq{diffw}). It expresses a standard relation between
Bernoulli polynomials.

Continuing this theme of rewriting, Adamchik's form of the GKBJ constants is,
  $$
  L_k=-\ze'(-k)-H_k\,\ze(-k)\,.
  \eql{GKBJ}
  $$
But I take this no further since, I will give, at another time, an extended analysis of this
approach and its generalisations as in [\pref{DandKi}], [\pref{KOW}] and elsewhere.

\section{\bf 9. Summary}

I have shown that Glaisher's numerics and Bendersky's recursive approach to the
generalised gamma functions were both anticipated by Jeffery whose work appeared a
couple of years after Kinkelin's paper, of which Jeffery seems unaware.

I have advocated the computation of the derivative of the Hurwitz \zf\ at the negative
integers via the generalised gamma function, rather then {\it vice versa}.

 \vglue .5 in
 \begin{ignore}
 \section{\bf Comparison with K and O}

I only consider $r=1$ and $\om=1$.
 Relation of notations.
   $$\eqalign{\wt\Ga_{1,k}(x+1)&=\log\Ga_k(x+1)\cr
   \log\Ga_{1,k}(x+1)&=\zeta_H'(-k,x+1)=\ze(-k,x+1)
   }
   \eql{nots}
   $$

Definition of $\wt\Ga$ is
  $$
\log\wt\Ga_{1,k}(x)=\log\Ga_{1,k}(x)-\ze'(-k)
\eql{kodef}
  $$
which is just (\peq{bas}).

K and O seem to make the definition
  $$
  \ze'_H(s,0)=\ze_R(s)=\ze(s)
  $$
because in Theorem 3 they have $\Ga_{1,k}(0)$ and give an example
  $$
\log\Ga_{1,1}(0,1)=\ze'(-1)
  $$
\cf (\peq{nots}).

Check that the recursion of Theorem 3 is the same as Bendersky's. Theorem 3 reads in K
and O's notation (setting $k\to k+1$,
  $$\eqalign{
(k+1)\int_0^xdt\,\bigg(&\log\Ga_{1,k}(t)-{1\over k+1}\ze(-k,t)\bigg)
=\log\Ga_{1,k+1}(x)-\log\Ga_{1,k+1}(0)
}
\eql{thm3}
  $$
This should agree with (\peq{recurs}), repeated here,
  $$
  \eqalign{
&(k+1)\int \La_k(x+1)\,dx+{1\over
k+1}\phi_{k+1}(x+1)+x\,H_k\,B_{k+1}=\La_{k+1}(x+1)\cr
&(k+1)\int \La_k(x+1)\,d(x+1)+{1\over
k+1}\phi_{k+1}(x+1)+x\,H_k\,B_{k+1}=\La_{k+1}(x+1)\cr
&(k+1)\int \La_k(t)\,dt+{1\over
k+1}\phi_{k+1}(t)+(t-1)\,H_k\,B_{k+1}=\La_{k+1}(t)\cr
&(k+1)\int_0^x \La_k(t)\,dt+{1\over
k+1}\phi_{k+1}(x)+(x-1)\,H_k\,B_{k+1}+C=\La_{k+1}(x)\cr
&C=\La_{k+1}(0)+H_kB_{k+1}=H_kB_{k+1}-L_k\cr
&(k+1)\int_0^x \La_k(t)\,dt+{1\over
k+1}\phi_{k+1}(x)+x\,H_k\,B_{k+1}=\La_{k+1}(x)+L_{k+1}\cr
}
\eql{recurs1}
  $$
or,
  $$\eqalign{
(k+1)\int_0^x \bigg(\log\Ga_k(t)+{H_k\,B_{k+1}\over k+1}-L_k\bigg)\,dt+{1\over
k+1}\phi_{k+1}(x)=\log\Ga_{k+1}(x)\cr
(k+1)\int_0^x \bigg(\log\Ga_k(t)+\ze'(-k)\bigg)\,dt+{1\over
k+1}\phi_{k+1}(x)=\log\Ga_{k+1}(x)\cr
(k+1)\int_0^x \ze'(-k,t)\,dt+{1\over
k+1}\phi_{k+1}(x)=\log\Ga_{k+1}(x)=\ze'(-k-1,x)-\ze'(-k-1,0),\cr
}
\eql{recurs2}
  $$

Now adjust the Bernoulli functions. Use
  $$
  {\phi_{k+1}(x)\over k+1}=\int_0^xdt \bigg(\phi_k(t)+{B_{k+1}\over k+1}\bigg)=
  \int_0^xdt {B_{k+1}(t)\over k+1}= \int_0^xdt \,\ze(-k,t)
  $$
to give, finally
  $$
  (k+1)\int_0^x \bigg(\ze'(-k,t)-{\ze(-k,t)\over k +1}\bigg)\,dt=\ze'(-k-1,x)-\ze'(-k-1,0),
  $$
which is (\peq{thm3}).

Correspondingly there is an easier proof of the recursion than Bendersky's, starting from
the standard relation (used in ({\pref{DoandKi}) \eg\ ),
  $$
  {\pa\over\pa w}\ze(s,w)=-s\ze(s+1,w)
  $$
so that
  $$
    \ze(s,w)-\ze(s,0)=-s\int_0^w dt \ze(s+1,t)
  $$
and
  $$
    \ze'(-k-1,w)-\ze'(-k-1,0)=-\int_0^w dt \ze(-k,t)+(k+1)\int_0^wdt\,\ze'(-k,t)
  $$

\end{ignore}

\begin{ignore}

To go beyond Jeffery's values, I look at the next case of $v_3(x)$ which is algebraically
more messy. Firstly, from (\peq{iter}) and (\peq{vee1}),
  $$\eqalign{
  v_3(1)=0&=6\int_0^1 dx\int_0^{x}dx'\,v_1(x')+{1\over2}-{3\over4\times6}+
  {3\over2}\varpi(2)+{1\over3}+\varpi(3)\cr
  &=6\int_0^1 dx\int_0^{x}dx'\int_0^{x'}dx''\,\bigg(\log\Ga(x''+1)+{x''(x''+1)\over2}-L_0x''\bigg)\cr
  &\hspace{*******}+{17\over24}+{3\over2}\varpi(2)+\varpi(3)\cr
  &=3\int_0^1 dx(1-x)^2\bigg(\log\Ga(x)+\log(x)+{x(x+1)\over2}-L_0x\bigg)\cr
  &\hspace{*******}+{17\over24}+{3\over2}\varpi(2)+\varpi(3)\cr
  &=3\int_0^1 dx(1-x)^2\log\Ga(x)+{L_0\over4}+{199\over120}
  +{17\over24}+{3\over2}\varpi(2)+\varpi(3)\cr
  &=3\int_0^1 dx(1-x)^2\log\Ga(x)+{L_0\over4}+{71\over30}+{3\over2}\varpi(2)+\varpi(3)\cr
  }
  \eql{vee3}
  $$
where I have immediately moved to $v_3(1)=0$. Next I again use (\peq{iprod2}) in the
integral
  $$
  3\int_0^1 dx(1-x)^2\log\Ga(x)=-\log w!+{3\over4}\log w -
  3\sum_{i=0}^{w-1}\int_0^1 dx\,(1-x)^2\log(x+i)\,.
  \eql{wint}
  $$
and, this time,
  $$\eqalign{
  3\int_0^1 dx\,(1-x)^2\log(x+i)=&i^3\log(i+1)
  +3i^2\log(i+1)+3i\log(i+1)\cr
  +&\log(i+1)-i^3\log i-3i^2\log i-3i\log i \cr&-i^2-{5i\over2}-{11\over6}
  }
  $$

Then
  $$\eqalign{
  &\sum_{i=0}^{w-1}\big(i^3\log(i+1)-i^3\log i\big)=
  \sum_{i=1}^{w}(3i-3i^2-1)\log i+w^3\log w\cr
   &\sum_{i=0}^{w-1}\big(3i^2\log (i+1)-3i^2\log i\big)=
    3\sum_{i=1}^{w}(1-2i)\log i+3w^2\log w\cr
    &\sum_{i=0}^{w-1}\big(3i\log (i+1)-3i\log i\big)=3w\log w\cr
  &\sum_{i=0}^{w-1}\bigg(-i^2-{5i\over2}-{11\over6}\bigg)=
  -{w^3\over3}-{7w^2\over4}-{5w\over12}
  }
  $$

Addition gives
  $$\eqalign{
  3\sum_{i=0}^{w-1}&\int_0^1 dx\,(1-x)^2\log(x+i)=\cr
  &\sum_{i=1}^{w-1}(-3i-3i^2+2)\log i+(w^3+3w^2+3w)\log w
  -{w^3\over3}-{7w^2\over4}-{5w\over12}
  }
  $$
and, going back to (\peq{wint}), I find for the integral needed in (\peq{vee3}),
  $$\eqalign{
  3\int_0^1 dx&(1-x)^2\log\Ga(x)\cr
  &=-\log w!+{3\over4}\log w-
  \sum_{i=1}^{w}(-3i-3i^2+2)\log i+(w^3+3w^2+3w)\log w\cr
  &=-3(-w+w\log w +{1\over2}\log w +L_0)\cr
&+{3\over4}\log w-
  \sum_{i=1}^{w}(-3i-3i^2+2)\log i+(w^3+3w^2+3w)\log w
  }
  \eql{wint2}
  $$

 I will start in some generality by defining a \zf,
  $$
   \ze(s,w)=\sum_\lambda{1\over(\lambda+w)^s}
  $$
where the $\lambda$   are eigenvalues of some differential (or pseudo--differential)
operator, $\caO$ and $w$ might be interpreted as a mass  or a Laplace transfer variable or
as just a handy parameter. In the special case when $\lambda=n$, ($n=1,2,\ldots)$
$\ze(s,w) $ is the Hurwitz \zf, $\ze_H(s,w)$.

The (asymptotic) expansion of $\ze(s,w)$ can easily be worked out, the coefficients being
related to those in the short--time expansion of the heat or cylinder--kernel, as the case
may be, of $\caO$.

In the Hurwitz case, these coefficients are proportional to Bernoulli numbers and the
expansion is essentially the old one of Mulholland for the partition function of a diatomic
molecule.

By differentiating this with respect to $s$, the expansion of $\ze'(-p,w)$ can be derived.
That for the important case of $p=0$ was given, \eg, in [\pref{DandC}], but the general
case is easily obtained. I will not write it out although it will be useful.

There is another route to this expansion which is essentially contained in [\pref{DandKi}].
where it is shown that,
  $$
  {\pa^n\over\pa w^n}\ze'(-n,w)=n!\,\big(\ze'(0,w)-H_n\,\ze(0,w)\big)
  $$
and, in the Hurwitz case,
  $$\eqalign
  {\pa^n\over\pa w^n}\ze_H'(-n,w)=n!\,
  \big(\log\big(\Ga(w)/\sqrt{2\pi}\big)-(1/2-w)H_n\big)
  $$

The idea now is to use the Stirling's formula on the right---hand side and integrate $n$
times. The right--hand side is
  $$
  n!\bigg(\big(w-1/2\big)\big(H_n+\log w\big)-w+
  \sum_{m=1}^\infty{B_{2m}\over(2m-1)\,2m\,w^{2m-1}}\bigg)\,.
  $$

  As an example take $n=1$.
\end{ignore}
 \vglue 20truept

 \noin{\bf References.} \vskip5truept
\begin{putreferences}
    \ref{EandM}{Espinosa,O. and Moll,V.H. {\it The Ramanujan Journal} {\bf6}
(2002) 449.}
   \ref{Gosper}{Gosper,R,W.}
     \ref{Jeffery3}{Jeffery, H.M. \qjm{4}{1861}{364}.}
     \ref{BarnesGf}{Barnes,E.W.{\it Quart.J.Pure and Applied Maths.}
     {\it 31} (1899) 264.}
   \ref{Jeffery2}{Jeffery, H.M. \qjm{5}{1862}{91}.}
   \ref{Jeffery}{Jeffery, H.M. \qjm{6}{1864}{82}.}
 \ref{Norlund}{N\"orlund,N.E. \am{43}{1922}{121}.}
 \ref{Glaisher}{Glaisher,J.W.L. {\it Messenger of Math.} {\bf 7} (1878) 43.}
 \ref{Dow}{Dowker,J.S. {\it Poweroids revisited -- an old symbolic approach.}ArXiv:1307.3150.}
 \ref{KOW}{Kurokawa,N,, Ochiai,H. and Wakayama, M. {\it J.Ramanujan Math.Soc.} {\bf21}
 (2006) 153.}
 \ref{KandO}{Kurokawa,N. and Ochiai,H. {\it Kodaira Math.J.} {\bf30} (2007) 195.}
    \ref{Elizalde}{Elizalde,E. {\it Math. of Comp.} {\bf 47} (1986) 347.}
    \ref{kink}{Kinkelin,H. {\it J.f.reine u. angew. Math. (Crelle)} {\bf 57} (1860)
   122.}
  \ref{holder}{O.H\"older {\it G\"ott. Nachrichten} (1886) 514-522.}
  \ref{alex}{W.P.Alexeiewsky {\it Leipzig Berichte} {\bf 46} (1894) 268-275.}
     \ref{bend}{Bendersky,L. \am{61}{1933}{263}.}
     \ref{Adamchik1}{Adamchik, V.S.
   {\it J. Comp. Appl. Math.} {\bf 100} (1998) 191.}
   \ref{Adamchik}{V.S.Adamchik, {\it Contributions to the theory of the
   Barnes function}, ArXiv: math. CA/0308086.}
    \ref{FandT}{Fradkin, E.S and Tseytin,A.A. \pl{B134}{1984}{301}.}
     \ref{CEZ2}{Cognola,G.,Elizalde,E. and Zerbini,S.  \cmp{237}{2003}{507}.}
   \ref{Tseytlin2}{Tseytlin, A. \np{877}{2013}{632}.}
   \ref{Tseytlin}{Tseytlin,A.A. \np{877}{2013}{598}.}
  \ref{Dowma}{Dowker,J.S. {\it Calculation of the multiplicative anomaly} ArXiv: 1412.0549.}
  \ref{CandH}{Camporesi,R. and Higuchi,A. {\it J.Geom. and Physics}
  {\bf 15} (1994) 57.}
  \ref{Allen}{Allen,B. \np{226}{1983}{228}.}
  \ref{Dowdgjms}{Dowker,J.S. \jpamt{48}{2015}{125401}.}
  \ref{Dowsphgjms}{Dowker,J.S. {\it Numerical evaluation of spherical GJMS determinants
  for even dimensions}, ArXiv:1310.0759.}
  \ref{CEZ}{Cognola,G.,Elizalde,E. and Zerbini,S.  \jpamt{48}{2015}{045203}.}
 \ref{CFM}{Castillo-Garate,V.,Friedman,E. and M\v{a}ntoiu,M. {\it The multiplicative
 anomaly of three or more commuting elliptic operstors}, ArXiv:1211.4117.}
  \ref{DowGJMS}{Dowker,J.S.  \jpa{44}{2011}{115402}.}
  \ref{Dowcmp}{Dowker,J.S. \cmp{162}{1994}{633}.}
\ref{Dowbfe}{Dowker,J.S. {\it The boundary F-theorem for free fields}, ArXiv:1407.5909.}
    \ref{DandKi}{Dowker,J.S. and Kirsten, K. {\it Comm. in Anal. and Geom.}
    {\bf7} (1999) 641.}
         \ref{Kassel}{Kassel,C. Seminaire Bourbaki  n. 708 (1988-89) 199.}
              \ref{Vardi}{Vardi,I. {\it SIAM J.Math.Anal.} {\bf 19} (1988) 493.}

  \ref{CandT2}{Copeland,E. and Toms,D.J. \cqg {3}{1986}{431}.}
   \ref{DoandKi} {Dowker.J.S. and Kirsten, K. {\it Analysis and Appl.}
       {\bf 3} (2005) 45.}
\begin{ignore}

  \ref{CandT}{Copeland,E. and Toms,D.J. \np {255}{1985}{201}.}
  \ref{Apps}{Apps,J.S. Thesis (University of Manchester, 1995).}
 \ref{Allen2}{Allen,B. PhD Thesis, University of Cambridge, 1984.}
 \ref{Chodos1}{Chodos,A. and Myers,E. \aop{156}{1984}{412}.}
     \ref{Guillarmou}{Guillarmou,C. \ajm{131}{2009}{1359}.}

     \ref{BaandDu}{Basar,G and Dunne,G.V. \jpa {43}{2010}{072002}.}
     \ref{AaandD}{Aros,R. and Diaz,D.E. {\it Determinant and Weyl anomaly of
     Dirac operator: a holographic derivation}, ArXiv:1111.1463.}
     \ref{BaandS}{B\"ar,C. and Schopka,S. The Dirac determinant of spherical
     space forms,\break {\it Geom.Anal. and Nonlinear PDEs} (Springer, Berlin, 2003).}
     \ref{QandC}{J.R.Quine and J.Choi, \rmjm {26}{1996}{719-729}.}
  \ref{QandC2}{J.R.Quine and J.Choi, Zeta regularized products and functional
  determinants on spheres \rmjm {26}{1996}{719-729}.}
  \ref{moller}{M{\o}ller,N.M. \ma {343}{2009}{35}.}
     \ref{KKY}{Kanemitsu,S., Kumagai,H. and Yoshimoto,M. {\it The Ramanujan
    J.} {\bf 5}(2001)5.}
     \ref{KandB}{Kamela,M. and Burgess,C.P. \cjp{77}{1999}{85}.}

     \ref{Dow20}{Dowker,J.S. \jmp{35}{1994}{6076}.}
     \ref{DowGJMSO}{Dowker,J.S. {\it Numerical evaluation of spherical GJMS operators}
     ArXiv: \break 1309.2873.}
     \ref{Adamchik1}{V.S.Adamchik, Polygamma functions of negative order
   {\it J. Comp. Appl. Math.} {\bf 100} (1998) 191-199.}
   \ref{Adamchik}{V.S.Adamchik, {\it Contributions to the theory of the
   Barnes function}, ArXiv: math. CA/0308086.}

    \ref{Onodera}{Onodera,K. \aim{224}{2010}{895}.}
    \ref{KKY}{Kanemitsu,S., Kumagai,H. and Yoshimoto,M. {\it The Ramanujan
    J.} {\bf 5}(2001)5.}
     \ref{Juhl}{Juhl,A. {\it On conformally covariant powers of the Laplacian}
     ArXiv: 0905.3992}
     \ref{Elphinstone}{Elphinstone,H.W. \qjm {2}{1858}{252}.}
       \ref{Branson}{Branson,T.P. \tams{347} {1995}{3671}.}
      \ref{Bromwich}{Bromwich, T.J.I'A. {\it Infinite Series},
  (Macmillan, London, 1926).}
   \ref{Loney}{Loney, S.L. {\it Plane Trigonometry} (CUP, Cambridge, 1893).}
    \ref{Hertzberg}{Hertzberg,M.P. \jpa{46}{2013}{015402}.}
     \ref{CaandW}{Callan,C.G. and Wilczek,F. \plb{333}{1994}{55}.}
    \ref{CaandH}{Casini,H. and Huerta,M. \plb{694}{2010}{167}.}
    \ref{Lindelof}{Lindel\"of,E. {\it Le Calcul des Residues} (Gauthier--Villars, Paris,1904).}
    \ref{CaandC}{Calabrese,P. and Cardy,J. {\it J.Stat.Phys.} {\bf 0406} (2004) 002.}
    \ref{MFS}{Metlitski,M.A., Fuertes,C.A. and Sachdev,S. \prB{80}{2009}{115122}.}
    \ref{Gromes}{Gromes, D. \mz{94}{1966}{110}.}
    \ref{Pockels}{Pockels, F. {\it \"Uber die Differentialgleichung $\De
  u+k^2u=0$} (Teubner, Leipzig. 1891).}
   \ref{Diaz}{Diaz,D.E. JHEP {\bf 0807} (2008) 103.}
    \ref{DandD}{Diaz,D.E. and Dorn,H. JHEP {\bf 0705} (2007) 46.}
  \ref{Minak}{Minakshisundaram,S. {\it J. Ind. Math. Soc.} {\bf 13} (1949) 41.}
    \ref{CaandWe}{Candelas,P. and Weinberg,S. \np{237}{1984}{397}.}
     \ref{Chodos1}{Chodos,A. and Myers,E. \aop{156}{1984}{412}.}
     \ref{ChandD}{Chang,P. and Dowker,J.S. \np{395}{1993}{407}.}
    \ref{LMS}{Lewkowycz,A., Myers,R.C. and Smolkin,M. {\it Observations on
    entanglement entropy in massive QFTs.} ArXiv:1210.6858.}
    \ref{Bierens}{Bierens de Haan,D. {\it Nouvelles tables d'int\'egrales
  d\'efinies}, (P.Engels, Leiden, 1867).}

    \ref{Dowren}{Dowker,J.S. \jpamt {46}{2013}{2254}.}
    \ref{Doweven}{Dowker,J.S. {\it Entanglement entropy on even spheres.}
    ArXiv:1009.3854.}
     \ref{Dowodd}{Dowker,J.S. {\it Entanglement entropy on odd spheres.}
     ArXiv:1012.1548.}
    \ref{DeWitt}{DeWitt,B.S. {\it Quantum gravity: the new synthesis} in
    {\it General Relativity} edited by S.W.Hawking and W.Israel (CUP,Cambridge,1979).}
    \ref{Nielsen}{Nielsen,N. {\it Handbuch der Theorie von Gammafunktion}
    (Teubner,Leipzig,1906).}
    \ref{KPSS}{Klebanov,I.R., Pufu,S.S., Sachdev,S. and Saddi,B.R.
    {\it JHEP} 1204 (2012) 074.}
    \ref{KPS2}{Klebanov,I.R., Pufu,S.S. and Safdi,B.R. {\it F-Theorem without
    Supersymmetry} 1105.4598.}
    \ref{KNPS}{Klebanov,I.R., Nishioka,T, Pufu,S.S. and Safdi,B.R. {\it Is Renormalized
     Entanglement Entropy Stationary at RG Fixed Points?} 1207.3360.}
    \ref{Stern}{Stern,W. \jram {79}{1875}{67}.}
    \ref{Gregory}{Gregory, D.F. {\it Examples of the processes of the Differential
    and Integral Calculus} 2nd. Edn (Deighton,Cambridge,1847).}
    \ref{MyandS}{Myers,R.C. and Sinha, A. \prD{82}{2010}{046006}.}
   \ref{RyandT}{Ryu,S. and Takayanagi,T. JHEP {\bf 0608}(2006)045.}

     \ref{Dowjmp}{Dowker,J.S. \jmp{35}{1994}{4989}.}
      \ref{Dowhyp}{Dowker,J.S. \jpa{43}{2010}{445402}.}
       \ref{HandW}{Hertzberg,M.P. and Wilczek,F. \prl{106}{2011}{050404}.}
      \ref{dowkerfp}{Dowker,J.S.\prD{50}{1994}{6369}.}
       \ref{Fursaev}{Fursaev,D.V. \plb{334}{1994}{53}.}
       \ref{Barnesa}{Barnes,E.W. {\it Trans. Camb. Phil. Soc.} {\bf 19} (1903) 374.}
  \ref{Barnesb}{Barnes,E.W. {\it Trans. Camb. Phil. Soc.}{\bf 19} (1903) 426.}
  \end{ignore}
\end{putreferences}

\bye